\documentclass[11pt]{article}

\usepackage{amsmath, amssymb, amscd}

\def\eqref#1{(\ref{#1})}
\newcommand{\goth}{\mathfrak}

\newcommand{\arrow}{{\:\longrightarrow\:}}
\newcommand{\Z}{{\Bbb Z}}
\newcommand{\C}{{\Bbb C}}
\newcommand{\R}{{\Bbb R}}

\newcommand{\6}{\partial}
\def\1{\sqrt{-1}\:}
\newcommand{\restrict}[1]{{\left|_{{\phantom{|}\!\!}_{#1}}\right.}}
\newcommand{\cntrct}                % contraction with a vector field
{\hspace{2pt}\raisebox{1pt}{\text{$\lrcorner$}}\hspace{2pt}}

\renewcommand{\c}[1]{{\cal #1}}
\newcommand{\calo}{{\cal O}}

\renewcommand{\tilde}{\widetilde}
\renewcommand{\bar}{\overline}
\renewcommand{\phi}{\varphi}
\renewcommand{\epsilon}{\varepsilon}
\renewcommand{\geq}{\geqslant}
\renewcommand{\leq}{\leqslant}

\newcommand{\End}{\operatorname{End}}
\newcommand{\Tot}{\operatorname{Tot}}
\newcommand{\Id}{\operatorname{Id}}

\newcommand{\Lie}{\operatorname{Lie}}

\newcommand{\Hol}{\operatorname{\it Hol}}

\newcommand{\Tr}{\operatorname{Tr}}

\newcommand{\comment}[1]{{}}

\def\blacksquare{\hbox{\vrule width 4pt height 4pt depth 0pt}}
\def\endproof{\blacksquare}

\makeatletter

\@ifundefined{Bbb}
     {\newcommand{\Bbb}[1]{{\mathbb #1}}}%
{}%     {\edef\Bbb#1{{\Bbb #1}}}

\newcommand{\ps@verbit}{%
  \renewcommand{\@oddhead}{%
          \scriptsize
          {Vanishing theorems for LCHK manifolds}
          \hfil\tiny {M. Verbitsky,  \ \ \ \ 15 February 2003}}
  \renewcommand{\@evenhead}{\@oddhead}
  \renewcommand{\@oddfoot}{\hfil\thepage\hfil}
  \renewcommand{\@evenfoot}{\@oddfoot}}
 
\newcounter{Mycounter}[section]
\newcounter{lemma}[section]
\renewcommand{\thelemma}{{Lemma \thesection.\arabic{lemma}}}
\newcommand{\lemma}{%
     \setcounter{lemma}{\value{Mycounter}}
     \refstepcounter{lemma}
     \stepcounter{Mycounter}
     {\bf \thelemma:\ }}

\newcounter{claim}[section]
\renewcommand{\theclaim}{{Claim \thesection.\arabic{claim}}}
\newcommand{\claim}{%
     \setcounter{claim}{\value{Mycounter}}
     \refstepcounter{claim}
     \stepcounter{Mycounter}
     {\bf \theclaim:\ }}

\newcounter{sublemma}[section]
\newcounter{corollary}[section]
\renewcommand{\thecorollary}{{Corollary \thesection.\arabic{corollary}}}
\newcommand{\corollary}{%
     \setcounter{corollary}{\value{Mycounter}}
     \refstepcounter{corollary}
     \stepcounter{Mycounter}
     {\bf \thecorollary:\ }}

\newcounter{theorem}[section]
\renewcommand{\thetheorem}{{Theorem \thesection.\arabic{theorem}}}
\newcommand{\theorem}{%
     \setcounter{theorem}{\value{Mycounter}}
     \refstepcounter{theorem}
     \stepcounter{Mycounter}
     {\bf \thetheorem:\ }}

\newcounter{conjecture}[section]
\newcounter{proposition}[section]
\renewcommand{\theproposition}
       {{Proposition \thesection.\arabic{proposition}}}
\newcommand{\proposition}{%
     \setcounter{proposition}{\value{Mycounter}}
     \refstepcounter{proposition}
     \stepcounter{Mycounter}
     {\bf \theproposition:\ }}

\newcounter{definition}[section]
\renewcommand{\thedefinition}
       {{Definition~\thesection.\arabic{definition}}}
\newcommand{\definition}{%
     \setcounter{definition}{\value{Mycounter}}
     \refstepcounter{definition}
     \stepcounter{Mycounter}
     {\bf \thedefinition:\ }}

\newcounter{example}[section]
\renewcommand{\theexample}{{Example \thesection.\arabic{example}}}
\newcommand{\example}{%
     \setcounter{example}{\value{Mycounter}}
     \refstepcounter{example}
     \stepcounter{Mycounter}
     {\bf \theexample:\ }}

\newcounter{remark}[section]
\renewcommand{\theremark}{{Remark \thesection.\arabic{remark}}}
\newcommand{\remark}{%
     \setcounter{remark}{\value{Mycounter}}
     \refstepcounter{remark}
     \stepcounter{Mycounter}
     {\bf \theremark:\ }}

\newcounter{problem}[section]
\newcounter{question}[section]
\@addtoreset{equation}{section}
\@addtoreset{footnote}{section}
\makeatother

\begin{document}

%%%%%%%%%%%%%%%%%%%%%%%%%%%%%%%%%%%%%%%%%%%%%%%%%%%%%%%%%%%%
\begin{center}
{\Large\bf
Vanishing theorems for locally
conformal \\[3mm] hyperk\"ahler manifolds
}
%%%%%%%%%%%%%%%%%%%%%%%%%%%%%%%%%%%%%%%%%%%%%%%%%%%%%%%%%%%%
\\[4mm]
Misha Verbitsky,\footnote{The author is 
partially supported by CRDF grant 
RM1-2354-MO02 and EPSRC grant  GR/R77773/01}
\\[4mm]
{\tt verbit@maths.gla.ac.uk, verbit@mccme.ru}
\end{center}

%%%%%%%%%%%%%%%%%%%%%%%%%%%%%%%%%%%%%%%%%%%%%%%%
{\small 
\hspace{0.15\linewidth}
\begin{minipage}[t]{0.7\linewidth}
{\bf Abstract} \\
Let $M$ be a compact locally
conformal hyperk\"ahler manifold.
We prove a version of Kodaira-Nakano
vanishing theorem for $M$. This is used
to show that $M$ admits no holomorphic differential
forms, and the cohomology of the structure sheaf
$H^i(\calo_M)$ vanishes for $i>1$. We also prove
that the first Betti number of $M$ is 1.
This leads to a structure theorem for locally
conformally hyperk\"ahler manifolds,
describing them in terms of 
3-Sasakian geometry.
Similar results are proven for
compact Einstein-Weyl locally
conformal K\"ahler manifolds.

\end{minipage}
}
%%%%%%%%%%%%%%%%%%%%%%%%%%%%%%%%%%%%%%%%%%%%%%%%

{
\small
\tableofcontents
}
%%%%%%%%%%%%%%%%%%%%%%%%%%%%%%%%%%%%%%%%%%%%%%%%

\section{Introduction}
\label{_Intro_Section_}

%%%%%%%%%%%%%%%%%%%%%%%%%%%%%%%%%%%%%%%%%%%%%%%%

The locally conformal K\"ahler manifolds and locally conformally
hyperk\"ahler manifolds were intensively studied throughout the
last 30 years; see \cite{_Dragomir_Ornea_} 
and \cite{_Ornea:LCHK_} for a survey of
known results and further reference. 
The key notion in this study is the notion
of a Vaisman manifold, also known
as generalized Hopf manifold (\ref{_Vaisman_Definition_}).
These manifolds were discovered by I. Vaisman and studied
in a big series of papers in early 1980-es (see 
\cite{_Vaisman:Dedicata_}, \cite{_Vaisman:Torino_},
and the bibliography in \cite{_Dragomir_Ornea_}).

We prove vanishing theorems for locally conformally
hyperk\"ahler manifolds and locally conformal K\"ahler
Einstein-Weyl manifolds and use these results to
obtain structure theorems. Similar vanishing theorems
were obtained in \cite{_Aleksandrov_Ivanov_}
using estimates on Ricci curvature.

In physics, 
the vanishing theorems for locally conformal K\"ahler manifolds
can be interpreted as conditions for existence of certain 
string compactifications (\cite{_Strominger_}, 
\cite{_Ivanov_Papado1_}, \cite{_Ivanov_Papado2_}).

%%%%%%%%%%%%%%%%%%%%%%%%%%%%%%%%%%%%%%%%%%%%%%%%%%%%%%%%%%%%
\subsection{Locally conformal hyperk\"ahler
manifolds, definition and examples}
\label{_LCHK_def_intro_Subsection_}
%%%%%%%%%%%%%%%%%%%%%%%%%%%%%%%%%%%%%%%%%%%%%%%%%%%%%%%%%%%%

Let $M$ be a smooth manifold equipped with operators $I, J, K \in \End (TM)$
satisfying the quaternion relations 
\[ I^2 = J^2=K^2 = -1, I\circ J = - J\circ I = K.
\]
Assume that 
the operators $I, J, K$ induce integrable complex structures on
$M$. Then $M$ is called {\bf hypercomplex}.
By a theorem 
of Obata (\cite{_Obata_}), a hypercomplex manifold admits
a unique torsion-free connection preserving $I, J, K$.
If the Obata connection preserves the metric $g$,
$(M,g)$ is called {\bf hyperk\"ahler}. Hyperk\"ahler
manifolds were introduced by E. Calabi
(\cite{_Calabi_}), and hypercomplex manifolds,
much later, by C.P.Boyer (\cite{_Boyer_}).

A hypercomplex manifold $M$ is called locally conformal hyperk\"ahler
(LCHK) if  the covering of $M$ is hyperk\"ahler, and the
monodromy transform preserves the conformal class of
a hyperk\"ahler metric. For a differently worded
definition, see Subsection \ref{_LCHK_Weyl_defi_Subsection_}.

The most elementary example of an LCHK manifold
is the Hopf manifold, defined as follows. Fix a quaternion
number $q\in {\Bbb H}$, $|q|>1$.
Consider the manifold $\tilde M = {\Bbb H}^n\backslash 0$,
and let $M:= \tilde M/\Z$, where $\Z$ acts on
$\tilde M$ by right quaternionic dilatations as 
\[ 
   (i, z) \arrow  z \cdot q^i, \ \ z\in {\Bbb H}^n, i\in \Z.
\]
This map is clearly compatible with the
hypercomplex structure and the conformal structure
on ${\Bbb H}^n$. Therefore, the manifold $M$ is locally
conformal hyperk\"ahler. 

More examples 
of LCHK manifolds are provided by the quaternionic
K\"ahler geometry (see e.g. \cite{_Besse:Einst_Manifo_}).
We shall not use these examples; the reader
not versed in quaternionic
K\"ahler geometry may skip the
next paragraph.

Given a quaternionic K\"ahler manifold $Q$ 
with positive scalar curvature (e.g.
the quaternionic projective space ${\Bbb H} P^n$)
one considers its Swann bundle $\c U(Q)$,
which fibers on $Q$ with a fiber $\C^2\backslash 0$
if $Q$ is Spin, or $\C^2\backslash 0/\{\pm 1\}$
if $Q$ is not Spin (\cite{_Swann_}). 
The total space of $\c U(Q)$
is hyperk\"ahler, and the natural dilatation 
map $\rho_t$ acts on $\c U(Q)$ preserving the
hypercomplex structure, and multiplies the metric
by a number. Take a quotient
$M:= \c U(Q)/\rho_{q^i}$, $i\in Z$, where $q>1$
is a fixed real number. If $Q$ is compact, then 
$M$ is also compact. By construction, $M$
is locally conformal hyperk\"ahler. 
When $Q={\Bbb H} P^n$, $\c U(Q) = \C^{2n}\backslash 0$,
and $M = \C^{2n}\backslash 0/\{ q^n\}$ is the 
Hopf manifold. 

%%%%%%%%%%%%%%%%%%%%%%%%%%%%%%%%%%%%%%%%%%%%%%%%%%%%%%%%%%%%
\subsection{Vanishing theorems for LCHK manifolds}
%%%%%%%%%%%%%%%%%%%%%%%%%%%%%%%%%%%%%%%%%%%%%%%%%%%%%%%%%%%%

In this paper we obtain several vanishing results for
LCHK (locally conformal hyperk\"ahler) manifolds
based on the same analytic arguments as the Kodaira-Nakano 
vanishing theorem. In particular, we obtain

\hfill

%%%%%%%%%%%%%%%%%%%%%%%%%%%%%%%%%%%%%%%%%%%%%%%%
\theorem\label{_Vanishing_intro_Theorem_}
Let $M$ be a compact LCHK manifold which is not hyperk\"ahler.
Consider $M$ as a complex manifold, with the complex structure
$I$ induced by the hypercomplex structure. Then 
\begin{description}
\item[(i)]
$M$ does not admit non-trivial global holomorphic forms
\item[(ii)] The cohomology of the structure sheaf of $M$
satisfy $H^i(\calo_M)=0$ for $i>1$, $\dim H^1(\calo_M)=1$.
\item[(iii)] The first Betti number of $M$ is 1.
\end{description}

{\bf Proof:} This is \ref{_holo_forms_vanish_Theorem_}, 
\ref{_cohomo_calo_M_Theorem_}, \ref{_h^1=1_Theorem_}. 
\endproof

\hfill

Using the Dolbeault spectral sequence, it is easy
to deduce \ref{_Vanishing_intro_Theorem_} (iii)  from
(i) and (ii). \ref{_Vanishing_intro_Theorem_} (iii) can
be proven directly by a simple geometric argument
(see Section \ref{_h^1=1_Appendix_}). 

\hfill

The same results are true 
for compact Einstein-Weyl 
locally conformally K\"ahler 
manifolds.

%%%%%%%%%%%%%%%%%%%%%%%%%%%%%%%%%%%%%%%%%%%%%%%%
\subsection{Geometry of LCHK manifolds}
%%%%%%%%%%%%%%%%%%%%%%%%%%%%%%%%%%%%%%%%%%%%%%%%

Consider a compact LCHK manifold $M$. Then
$M$ admits a special metric, discovered by P. Gauduchon
(\ref{_Gaud_metric_Definition_}). The Gauduchon
metric is defined as follows.

Let $\tilde M$ be the hyperk\"ahler covering of $M$.
Since the deck transform acs on $\tilde M$ preserving
the conformal class of the metric, the manifold $M$
is equipped with the canonical conformal structure
$[g]$. The Obata connection $\nabla$ on a hyperkaehler
manifold coincides with the Levi-Civita connection.
Therefore, $\nabla$ preserves
the conformal class $[g]$. We obtain that a parallel
transport along $\nabla$ multiplies a metric
$g\in [g]$ by a number. Therefore,
$\nabla(g) = g\otimes \theta$, where
$\theta$ is a 1-form, called
{\bf the Lee form of $(M,g)$.}

A metric $g\in [g]$ is called {\bf a Gauduchon metric}
if $\theta$ is co-closed. Gauduchon proved that this
metric exists and is unique up to a constant
multiplier if we fix the connection and the conformal
class (\ref{_Gaud_metric_exists_Lemma_}). 

If $M$ is an LCHK manifold equipped with a Gauduchon metric,
then the Lee form $\theta$ is parallel 
with respect to the Levi-Civita connection
(\ref{_LCHK_Vaisman_Theorem_}).  
The corresponding vector field
$\theta^\sharp$ is obviously Killing.
Moreover, the flow, associated with 
$\theta^\sharp$, is compatible with the hypercomplex structure
on $M$ (\ref{_flow_along_theta_automo_Proposition_}).
If we lift $\theta^\sharp$ to a hyperk\"ahler
covering $\tilde M$, the corresponding flow
multiplies the hyperk\"ahler metric by constant.
This allows one to construct the K\"ahler
potential on $\tilde M$ explicitly in terms
of the Lee form (\ref{_kahler_pote_on_Vaism_Proposition_}). 

Applying \ref{_kahler_pote_on_Vaism_Proposition_}
to different induced complex structures, we find
that $\tilde M$ is equipped with {\bf a hyperk\"ahler
potential}, that is, a function which serves
as a K\"ahler potential for all induced 
complex structures.

Hyperk\"ahler manifolds admitting hyperk\"ahler
potential were studied by A. Swann \cite{_Swann_}.
These manifolds are deeply related to quaternionic 
K\"ahler geometry. Take a 4-dimensional foliation
$\Phi$ generated by the gradient $\theta^\sharp$ of the
hyperk\"ahler potential and 
$I (\theta^\sharp), J (\theta^\sharp), K (\theta^\sharp) $.
This foliation is integrable, flat and completely geodesic.
The leaf space of $\Phi$ is quaternionic K\"ahler.
This leads to nice structure
theorems for hyperk\"ahler manifolds 
admitting hyperk\"ahler potential (see
\cite{_Swann_}).

The LCHK geometry is much more delicate, due
to possible global irregularities of the foliation
$\Phi$. 

Let $M$ be an LCHK manifold equipped with a Gauducon metric,
$\theta$ is Lee form, and $\Phi$ the 4-dimensional foliation
generated by $\theta^\sharp$ ,
$I (\theta^\sharp), J (\theta^\sharp), K (\theta^\sharp) $.
One can speak of the leaf space of $\Phi$ if 
every point $x\in M$ has a 
neighbourhood $U\subset M$ 
such that every leaf of  $\Phi$
meets $U$ in finitely many connected components.
In this case the foliation $\Phi$ and the manifold
$M$ is called {\bf quasiregular}. The leaf
space $Y$ of a quasiregular foliation is an orbifold.

\ref{_Vanishing_intro_Theorem_} (iii)
(the equality $h^1(M)=1$) is proven for 
quasiregular LCHK manifolds
(\cite{_Ornea_Piccini_}).

Given a compact quasiregular LCHK manifold, the leaf space $Q$ of 
$\Phi$ is a quaternionic K\"ahler orbifold, and $M$ is fibered
over $Q$ with fibers which are isomorpic to Hopf surfaces
(\cite{_Ornea:LCHK_}, \cite{_Ornea_Piccini_}). 

In this paper we give a similar structure theorem
for LCHK manifolds with no quasiregularity assumption.

%%%%%%%%%%%%%%%%%%%%%%%%%%%%%%%%%%%%%%%%%%%%%%%%%%%%%%%%%%%%
\subsection{3-Sasakian geometry and structure theorem for LCHK manifolds}
%%%%%%%%%%%%%%%%%%%%%%%%%%%%%%%%%%%%%%%%%%%%%%%%%%%%%%%%%%%%

It is more convenient to speak of 3-Sasakian manifolds than
of quaternionic K\"ahler orbifolds. A 3-Sasakian manifold
is locally a bundle over a quaternionic K\"ahler orbifold,
with a fiber isomorphic to $SU(2)/\Gamma$, where $\Gamma\subset SU(2)$
is a finite group. 

One defines 3-Sasakian manifolds as follows.

Given a Riemannian manifold $(X, g)$,
the {\bf cone} $C(X)$ of $X$ is defined as
a Riemannian manifold $X \times \R^{>0}$
with the metric $t^2 g + dt^2$, where $t$ is the 
parameter in $\R^{>0}$.

A {\bf 3-Sasakian structure} on $X$ is a hyperk\"ahler
structure on $C(X)$, defined in such a way
the map $(x, t) \arrow x, qt$ is holomorphic,
for all $q \in \R^{>0}$. 3-Sasakian manifolds were
discovered in late 1960-ies (see \cite{_Udriste_} and \cite{_Kuo_}),
and studied extensively in mid-1990-ies
by Boyer, Galicki and Mann (\cite{_Boyer_Galicki_Mann_});
see the excellent survey \cite{_Boyer_Galicki_}.

The 3-Sasakian (and, more generally,
Einstein-Sasakian) manifolds can 
also be obtained as circle bundles over
Einstein Fano orbifolds. This gives a way to construct
extensive lists of examples of Einstein-Sasakian manifolds 
(see e.g.  \cite{_Boyer_Galicki:new-Einstein_} 
\cite{_Boyer_Galicki_Nakamaye2_}).

Given a compact LCHK manifold
$M$, its universal covering $\tilde M$
is hyperk\"ahler. Using the explicit description of the
hyperk\"ahler metric in terms of the Gauduchon 
metric (\ref{_kahler_pote_on_Vaism_Proposition_}), 
we find that $\tilde M$ is a cone manifold:
$\tilde M = C(X)$, where $X$ is 3-Sasakian
(\ref{_Vaisman_decompo_Sasa_Proposition_}). 

Fix $q\in \R^{>1}$. 
Consider the equivalence relation on $C(X)$
generated by $(x, t) \sim (x, qt)$.
The quotient $C(X)/\sim_q$ is clearly 
an LCHK manifold.

The structure theorem \ref{_LCHK_Structure_Theorem_}
describes any compact LCHK manifold in terms of 3-Sasakian
manifolds, as follows. We show that $M \cong C(X)/\sim_{\phi, q}$,
where $q\in \R^{>1}$, $\phi:\; X \arrow X$ is a 3-Sasakian isometry,
and $\sim_{\phi, q}$ an equivalence relation on $C(X)$
generated by $(x, t) \sim (\phi(x), qt)$.

%%%%%%%%%%%%%%%%%%%%%%%%%%%%%%%%%%%%%%%%%%%%%%%%
\subsection{Subvarieties of Vaisman manifolds}
%%%%%%%%%%%%%%%%%%%%%%%%%%%%%%%%%%%%%%%%%%%%%%%%

We also obtain the following application.

\hfill

%%%%%%%%%%%%%%%%%%%%%%%%%%%%%%%%%%%%%%%%%%%%%%%%

%%%%%%%%%%%%%%%%%%%%%%%%%%%%%%%%%%%%%%%%%%%%%%%%
\proposition\label{_subva_intro_Proposition_}
Let $M$ be a compact Vaisman manifold,
and $X\subset M$ a closed complex subvariety.
Then $X$ is tangent to the canonical foliation $\Xi$
in all its smooth points. In particular, if $X$ is smooth,
then $X$ is also a Vaisman manifold.

\hfill

{\bf Proof:} See \ref{_subva_Vaisman_Proposition_}. \endproof

\hfill

For smooth manifolds, \ref{_subva_intro_Proposition_}
is proven by K. Tsukada, see \cite{_Tsukada:subva_}.

%%%%%%%%%%%%%%%%%%%%%%%%%%%%%%%%%%%%%%%%%%%%%%%%

\section{Locally
conformal hyperk\"ahler manifolds}
\label{_LCHK_Section_}

%%%%%%%%%%%%%%%%%%%%%%%%%%%%%%%%%%%%%%%%%%%%%%%%

In this Section we give the definitions
and site some results related to Weyl geometry and
locally conformal hyperk\"ahler manifolds. We follow
\cite{_Ornea:LCHK_}.

%%%%%%%%%%%%%%%%%%%%%%%%%%%%%%%%%%%%%%%%%%%%%%%%
\subsection{Weyl structures}
\label{_Weyl_Subsection_}
%%%%%%%%%%%%%%%%%%%%%%%%%%%%%%%%%%%%%%%%%%%%%%%%

For an introduction to Weyl geometry and further
reference on the subject see e.g. 
\cite{_Calderbank_Pedersen_}.

\hfill

%%%%%%%%%%%%%%%%%%%%%%%%%%%%%%%%%%%%%%%%%%%%%%%%
\definition
Let $(M, g)$ be a Riemannian manifold, and $\nabla$ a torsion-free 
connection on $M$. Assume that $\nabla$ preserves the conformal
class of $g$, that is,
\begin{equation}\label{_Weyl_conne_formula_Equation_}
   \nabla g = g \otimes \theta \in S^2 T^* M \otimes T^* M,
\end{equation}
where $\theta$ is a 1-form. Then $(M, g, \nabla, \theta)$
is called {\bf a Weyl manifold}. If $d \theta=0$,
$M$ is called {\bf a closed Weyl manifold}. 
The form $\theta$ is called {\bf the Lee form},
or {\bf the Higgs field of $M$}. A torsion-free connection
which satisfies \eqref{_Weyl_conne_formula_Equation_}
is called {\bf Weyl connection}.

\hfill

%%%%%%%%%%%%%%%%%%%%%%%%%%%%%%%%%%%%%%%%%%%%%%%%
\remark
Let $M$ be a Riemannian manifold equipped
with a Levi-Civita connection. Then $M$ is a 
Weyl manifold, with trivial $\theta$.

\hfill

Given a Weyl manifold $(M, g, \nabla, \theta)$
and a function $\zeta:\; M \arrow \R$, we observe that
\[ 
  (M, e^\zeta g, \nabla, \theta+d\zeta)
\]
is also a Weyl manifold. Indeed,  
\[ \nabla(e^\zeta g) = e^\zeta \nabla(g) + e^\zeta g\otimes d\zeta =
   e^\zeta g\otimes (\theta+d\zeta).
\]

\hfill

%%%%%%%%%%%%%%%%%%%%%%%%%%%%%%%%%%%%%%%%%%%%%%%%
\definition\label{_confo_equi_Definition_}
In the above assumptions, 
the Weyl manifolds \[ (M, g, \nabla, \theta) \text{\ \ and \ \ }
(M, e^\zeta g, \nabla, \theta+d\zeta)\] are called
{\bf globally conformal equivalent}. 

\hfill

Let $(M, g, \nabla, \theta)$ be a Weyl manifold,
and $L_\R$ a trivial 1-dimensional real bundle on $M$. Denote
the trivial connection on $L$ by $\nabla_{tr}$. Consider
a connection 
${}^L\nabla := \nabla_{tr} - \frac{\theta}{2}$ 
on $L$.  If $(M, g, \nabla, \theta)$ is closed,
the bundle $(L, {}^L\nabla)$ is flat, because
$({}^L\nabla)^2 = d\theta=0$. Denote by $L$ the
complexification of $L_\R$, with the induced connection.

\hfill 

%%%%%%%%%%%%%%%%%%%%%%%%%%%%%%%%%%%%%%%%%%%%%%%%
\definition \label{_weight_bu_Definition_}
The bundle $(L, {}^L\nabla)$ is called
{\bf the weight bundle} of a Weyl manifold $M$.
The weight bundle is equipped with a trivialization
$\lambda$.

\hfill

\hfill

%%%%%%%%%%%%%%%%%%%%%%%%%%%%%%%%%%%%%%%%%%%%%%%%
\claim\label{_weight_confo_equi_Claim_}
Let $(M, g, \nabla, \theta)$ and
$(M,  g', \nabla, \theta')$
be conformal equivalent closed Weyl manifolds,
and $L, L'$ the corresponding weight bundles. Then 
$L, L'$ are isomorphic as flat bundles.

\hfill

{\bf Proof:} Write $g' = e^\zeta g$, $\theta' = \theta + d\zeta$.
Let $\lambda$, $\lambda'$ be the sections of $L$, $L'$ inducing the standard
trivialization. Then $e^\zeta\lambda$ satisfies
\[ 
  {}^L \nabla(e^\zeta\lambda) = \theta + d\zeta(e^\zeta\lambda) =
  \theta' + {}^L \nabla(e^\zeta\lambda).
\]
\endproof

\hfill

%%%%%%%%%%%%%%%%%%%%%%%%%%%%%%%%%%%%%%%%%%%%%%%%
\remark
Using the trivialization of $L$, 
we may consider $g$ as a 2-form on $M$ with 
values in $L^{\otimes 2}$.
Then $g$ is a parallel section of 
$S^2 T^*M \otimes L^{\otimes 2}$.

\hfill

%%%%%%%%%%%%%%%%%%%%%%%%%%%%%%%%%%%%%%%%%%%%%%%%
\remark \label{_weight_determi_Remark_}
Let $n=\dim_\R M$. Denote by $\c P$ the 
principal $GL(n)$-bundle $\c P$ associated with $TM$.
Clearly, $L$ is a line bundle associated with 
the representation $(\det T^* M)^{-\frac{1}{n}}$
of $\c P$.

\hfill

Let  $(M, g, \nabla, \theta)$ be a closed Weyl manifold,
and $L$ its weight bundle. The natural connection ${}^L\nabla$
in $L$ is flat because $M$ is closed. Consider a covering 
$(\tilde M, \tilde g, \tilde \nabla, \tilde \theta)$
of $M$, such that the lift $\tilde L$ of $L$ to $\tilde M$ has trivial monodromy.
Let $\zeta_0$ be a ${}^L\nabla$-parallel section of $\tilde L$, 
$\zeta_0\neq 0$, and $\lambda$ the trivialization defined
on $L$ as on a weight bundle. The quotient
$\zeta:=\frac{\zeta_0}{\lambda}$ satisfies
\[ 
   d \zeta = \frac{\nabla_{tr}(\zeta_0)}{\lambda} = \frac \theta 2
\]
because $\nabla_{tr} - {}^L\nabla = \frac \theta 2$.
Therefore, the manifolds
$(M, g, \nabla, \theta)$ and \linebreak 
$(M, e^{2\zeta} g, \nabla, 0)$ are conformally equivalent.
We obtain the following claim.

\hfill

%%%%%%%%%%%%%%%%%%%%%%%%%%%%%%%%%%%%%%%%%%%%%%%%
\claim \label{_theta_exact_confo_Riemann_Claim_}
Let  $(M, g, \nabla, \theta)$ be a closed Weyl manifold,
and $(L, {}^L\nabla)$ its weight bundle. 
Assume that the monodromy of $L$ is trivial. Then 
$(M, g, \nabla, \theta)$ is conformal equivalent
to a Riemannian manifold equipped with a Levi-Civita connection.
Conversely, if $(M, g, \nabla, \theta)$ is conformal equivalent
to a Riemannian manifold, then its weight bundle 
$L$ has trivial monodromy. 

\endproof

%%%%%%%%%%%%%%%%%%%%%%%%%%%%%%%%%%%%%%%%%%%%%%%%
\subsection{Locally conformal hyperk\"ahler (LCHK) manifolds: the definition}
\label{_LCHK_Weyl_defi_Subsection_}
%%%%%%%%%%%%%%%%%%%%%%%%%%%%%%%%%%%%%%%%%%%%%%%%

%%%%%%%%%%%%%%%%%%%%%%%%%%%%%%%%%%%%%%%%%%%%%%%%
\definition
Let $M$ be a hypercomplex manifold, and $\nabla$ the 
Obata connection on $M$ (see Subsection 
\ref{_LCHK_def_intro_Subsection_}). Assume that
$M$ admits a quaternionic Hermitian metric 
$g$ and a closed form $\theta$, such that 
$(M, g, \nabla, \theta)$ is a closed Weyl 
manifold. Then $M$ is called {\bf 
locally conformal hyperk\"ahler manifold},
or {\bf LCHK manifold}. 

\hfill

Let $M$ be an LCHK manifold, $L$ the corresponding
weight bundle, and $(\tilde M, \tilde g, \tilde\nabla)$ 
the covering associated with the monodromy 
representation of $L$. 
In \ref{_theta_exact_confo_Riemann_Claim_} we constructed
a Riemannian metric $g' = e^{2\zeta} \tilde g$ which is
preserved by $\tilde \nabla$. 
Then $\tilde \nabla$ is the Levi-Civita connection for
$(\tilde M, g')$. On the other hand, $\tilde \nabla$
preserves the quaternion action. As we have mentioned
in the Introduction (Subsection 
\ref{_LCHK_def_intro_Subsection_}), this implies that,
the manifold $(\tilde M, g')$ is hyperk\"ahler.
We obtain

\hfill

%%%%%%%%%%%%%%%%%%%%%%%%%%%%%%%%%%%%%%%%%%%%%%%%
\claim \label{_LCHK_w_trivial_monodro_Claim_}
Let $M$ be an LCHK (locally conformal hyperk\"ahler) manifold,
$L$ its weight bundle, and $\tilde M$ the covering of $M$ 
associated with its monodromy. Then $\tilde M$
is equipped with a hyperk\"ahler metrics, which 
is determined uniquely up to a constant multiplier.

\endproof

\hfill

The converse statement is also true.

\hfill

%%%%%%%%%%%%%%%%%%%%%%%%%%%%%%%%%%%%%%%%%%%%%%%%
\proposition \label{_LCHK_as_hc_w_hk_cover_Proposition_}
Let $M$ be a hypercomplex manifold, and $\tilde M$ its universal
covering. Assume that $\tilde M$ is equipped with a hyperk\"ahler
metric $g'$ in such a way that for any $\gamma \in \pi_1(M)$ the
corresponding deck transform $k_\gamma:\; \tilde M\arrow \tilde M$
multiplies $g'$ by a scalar $c(\gamma)\in \R^{>0}$.
Then $M$ admits a LCHK metric, which is determined
uniquely up to conformal equivalence. 

\hfill

{\bf Proof:} The map $\gamma \arrow c(\gamma)$
defines a 1-dimensional representation 
$c:\; \pi_1(M) \arrow \R^{>0}$.
Let $L$ be the corresponding flat bundle on $M$. Since
$\R^{>0}$ is contractible, $L$ is topologically trivial. 
Pick a nowhere degenerate section $\lambda$ of $L$.
The hyperk\"ahler metrics $\tilde g$ on $\tilde M$
can be considered as a map $S^2 T^* M \arrow L^{\otimes 2}$.
Using $\lambda^2$ as a trivialization of $L^{\otimes 2}$, we obtain a 
hypercomplex Hermitian metrics $g_\lambda$ on $M$.
Clearly, 
\[ 
  \nabla g_\lambda = \nabla(\lambda^2 \tilde g)= 
  \lambda^2 \tilde g \otimes \nabla(\lambda^2) = 
  g_\lambda \otimes (- 2 \theta_L),
\]
where $\theta_L$ is the connection form of $L$
associated with the trivialization $\lambda$. 
Therefore, $(M, g_\lambda, \nabla, - 2 \theta_L)$
is an LCHK manifold. 
\endproof

\hfill

%%%%%%%%%%%%%%%%%%%%%%%%%%%%%%%%%%%%%%%%%%%%%%%%
\remark
\ref{_LCHK_as_hc_w_hk_cover_Proposition_}
gives a nice interpretation of LCHK geometry, which
is much more clear than the usual approach. We shall sometimes
implicitly use \ref{_LCHK_as_hc_w_hk_cover_Proposition_} instead
of the definition.

%%%%%%%%%%%%%%%%%%%%%%%%%%%%%%%%%%%%%%%%%%%%%%%%

\section{Vaisman manifolds}
\label{_Vaisman_Section_}
%%%%%%%%%%%%%%%%%%%%%%%%%%%%%%%%%%%%%%%%%%%%%%%%

In this Section we present some introductory material in 
locally conformal K\"ahler geometry. For 
more details and a bibliography the reader is 
referred to \cite{_Dragomir_Ornea_}.

\hfill

%%%%%%%%%%%%%%%%%%%%%%%%%%%%%%%%%%%%%%%%%%%%%%%%
\definition
Let $M$ be a complex manifold equipped with a Hermitian metric $\omega$
and a closed Weyl connection $\nabla$, 
\[ \nabla(\omega) = \omega \otimes \theta, \ \ \theta\in \Lambda^1(M). \]
Assume that $\nabla$ preserves the complex structure operator:
$\nabla (I) =0$.
Then $M$ is called {\bf locally conformal K\"ahler (LCK)
manifold}, and $\theta$ is called {\bf the Lee form} of $M$.

\hfill

%%%%%%%%%%%%%%%%%%%%%%%%%%%%%%%%%%%%%%%%%%%%%%%%
\remark
Given an LCHK (locally conformal hyperk\"ahler)
manifold $M$, we obtain an LCK structure $(M,J)$
for every quaternion $J, J^2=-1$. 

\hfill

The basic working tool of locally conformal K\"ahler geometry
is the following lemma of P. Gauduchon (\cite{_Gauduchon_1984_}).

\hfill

%%%%%%%%%%%%%%%%%%%%%%%%%%%%%%%%%%%%%%%%%%%%%%%%
\lemma\label{_Gaud_metric_exists_Lemma_}
\cite{_Gauduchon_1984_}
Let $M$ be a compact, oriented conformal manifold, \linebreak $\dim_\R M >2$.
For any Weyl connection preserving the conformal structure, there
exists a metric $g_0$ in this conformal class, such that
the corresponding Lee form $\theta$ is co-closed 
with respect to $g_0$. This metric is
unique up to a constant multiplier.

 \endproof

\hfill

%%%%%%%%%%%%%%%%%%%%%%%%%%%%%%%%%%%%%%%%%%%%%%%%
\definition\label{_Gaud_metric_Definition_}
Let $M$ be a conformal manifold equipped with
a Weyl connection preserving the conformal structure $[g]$,
and $g_0$ the metric in the conformal class $[g]$, such that
the corresponding Lee form $\theta$ is co-closed.
Then $g_0$ is called {\bf a Gauduchon metric}.

\hfill

%%%%%%%%%%%%%%%%%%%%%%%%%%%%%%%%%%%%%%%%%%%%%%%%
\remark
If $(M, g, \nabla, \theta)$ is a closed Weyl manifold, 
and $g$ is a Gauduchon metric, then $\theta$ is harmonic.

\hfill

For compact LCHK manifolds, the Gauduchon metric
has parallel Lee form, in the following sense.

\hfill

%%%%%%%%%%%%%%%%%%%%%%%%%%%%%%%%%%%%%%%%%%%%%%%%
\theorem\label{_LCHK_Vaisman_Theorem_}
\cite{_Pedersen_Poon_Swann:Einstein_Weyl_}
Let $M$ be an LCHK manifold equipped with a Gauduchon
metric $g$, and $\theta$ its Lee form. Then
$\nabla_g(\theta)=0$, where $\nabla_g$
is the Levi-Civita connection associated with $g$.

\endproof

\hfill

%%%%%%%%%%%%%%%%%%%%%%%%%%%%%%%%%%%%%%%%%%%%%%%%
\definition\label{_Vaisman_Definition_}
Let $(M,g, \nabla, \theta)$ be an LCK manifold, and $\nabla_g$
its Levi-Civita connection. We say that $M$ is an {\bf LCK
manifold with parallel Lee form}, or {\bf Vaisman manifold},
if $\nabla_g (\theta) =0$. If $\theta\neq 0$, then
after rescaling, we may always 
assume that $|\theta|=1$. 
Unless otherwise stated,
we shall assume implicitly that $|\theta|=1$
for all Vaisman manifolds we consider. 

\hfill

Vaisman manifolds were introduced by I. Vaisman
under the name ``generalized Hopf manifolds''
in a big series of papers (see e.g. \cite{_Vaisman:Dedicata_}, 
\cite{_Vaisman:Torino_}) and studied 
extensively since then.

\hfill

%%%%%%%%%%%%%%%%%%%%%%%%%%%%%%%%%%%%%%%%%%%%%%%%
\remark
If $\nabla_g \theta=0$, then $\theta$ is harmonic with respect
to $g$. Therefore, the metric $g$ is automatically a Gauduchon
metric.

\hfill

%%%%%%%%%%%%%%%%%%%%%%%%%%%%%%%%%%%%%%%%%%%%%%%%
\example \label{_Hopf_Example_}
Fix a quaternion $q\in {\Bbb H}$, $|q|>1$.
Let $M := {\Bbb H}^n \backslash 0 / \sim_q$,
where $\sim_q$ is an equivalence relation generated
by $z \sim_q qz$. This manifold is called 
{\bf a Hopf manifold}. Since the multiplication
by $q$ preserves the flat connection
$\nabla_{\text{fl}}$, this connection can
be obtained as a pullback of a connection $\nabla$ on $M$.
Since $\nabla_{\text{fl}}$ preserves the conformal class
of a flat metric, $(M, \nabla)$ is a Weyl manifold.
On the other hand, $M$ is by construction
a hypercomplex manifold, and its covering 
a hyperk\"ahler one. By \ref{_LCHK_as_hc_w_hk_cover_Proposition_}
this implies that $M$ is an LCHK (locally conformal hyperk\"ahler)
manifold. Topologically, we have $M = S^{4n-1}\times S^1$.
The Gauduchon metric is the standard metric on 
$S^{4n-1}\times S^1$, and the form $\theta$
is the coordinate form lifted from $S^1$. Therefore,
$\theta$ is parallel.

\hfill

For other examples of Vaisman manifolds, see
e.g. \cite{_Gauduchon_Ornea_}, \cite{_Belgun_}, 
\cite{_Kamishima_Ornea_}.

%%%%%%%%%%%%%%%%%%%%%%%%%%%%%%%%%%%%%%%%%%%%%%%%%%%%%%%%%%%%
\section{K\"ahler potential on Vaisman manifolds}
%%%%%%%%%%%%%%%%%%%%%%%%%%%%%%%%%%%%%%%%%%%%%%%%%%%%%%%%%%%%

We construct a K\"ahler potential on a Vaisman manifold with exact
Lee form (see \cite{_Vaisman:Torino_}). 

\hfill

%%%%%%%%%%%%%%%%%%%%%%%%%%%%%%%%%%%%%%%%%%%%%%%%%%%%%%%%%%%%%%%%%%%%%%%%
\proposition\label{_flow_along_theta_automo_Proposition_}
Let $M$ be an LCK manifold with
parallel Lee form $\theta$, and $\theta^\sharp$
the vector field dual to $\theta$. Consider a diffeomorphism
flow $\psi_t$ associated with $\theta^\sharp$.
Then $\psi_t$ acts on $M$ preserving the LCK structure.

\hfill

{\bf Proof:} For a more detailed proof 
see e.g. \cite{_Dragomir_Ornea_}.

Since $\theta$ is parallel, $\theta^\sharp$ is a parallel
vector field. Therefore, $\theta^\sharp$ is Killing, and
$\psi_t$ acts on $M$ by isometries.

On the other hand, $\psi_t$ is a geodesic flow along $\theta^\sharp$,
therefore its differential $d \psi_t:\; T_x X \arrow T_{\psi_t(x)} X$
is equal to the parallel transport along the geodesics associated
with $\theta^\sharp$. Since the holonomy of $\nabla_g$ is
contained in $U(n) \cdot \R\subset GL(n, \C)$,
we obtain that $d\psi_t$ is $\C$-linear. Therefore,
$\psi_t$ is holomorphic. We find that $\psi_t$ preserves
the complex and the Hermitian structure on $M$. The
Weyl connection $\nabla$ can be written explicitly in terms of
the Levi-Civita connection and the Lee form (see e.g. \cite{_Ornea:LCHK_}, 
Definition 1.1)
\begin{equation}\label{_nabla_Weyl_via_nabla_g_Equation_}
  \nabla = \nabla_g - \frac{1}{2}(\theta\otimes Id + Id \otimes\theta - g \otimes \theta^\sharp),
\end{equation}
where $\nabla_g$ is the Levi-Civita connection on $M$.
Since $\psi_t$ preserves $\theta$, $g$ and $\nabla_g$, 
we obtain that $\psi_t$ preserves $\nabla$. 
\ref{_flow_along_theta_automo_Proposition_} is proven.
\endproof

\hfill

The Lee form $\theta$ is by definition closed. Passing to
a covering if necessary, we may assume that it is exact:
$\theta = dt$. Write $r = e^{-t}$. In the \ref{_Hopf_Example_},
$r$ is the radius function. In this case, $r$ is obviously
a K\"ahler potential of $M$.

\hfill

%%%%%%%%%%%%%%%%%%%%%%%%%%%%%%%%%%%%%%%%%%%%%%%%
\definition
Let $M$ be an LCK manifold with exact Lee form $\theta = dt$.
The function $r:= e^{-t}$ is called {\bf the potential} of $M$. 
Clearly, $r$ is defined uniquely, up to a positive constant
multiplier.

\hfill

%%%%%%%%%%%%%%%%%%%%%%%%%%%%%%%%%%%%%%%%%%%%%%%%
\claim\label{_kah_form_on_LCK_w_pote_Claim_}
Let $(M, g, \nabla)$ be an LCK manifold with 
exact Lee form $\theta$, $r$ its potential
and $\omega\in \Lambda^{1,1}(M)$ the Hermitian
form of $(M, g)$. Then $r\omega$ is a K\"ahler 
form.

\hfill

{\bf Proof:} Clearly, 
\[ 
  d (r\omega) = - \theta \wedge r\omega + r \cdot d\omega =
              -\theta \wedge r \omega + r \theta \wedge \omega =0
\]
On the other hand, $r\omega$ is positive definite, because
$r$ is a positive function. A closed positive definite
$(1,1)$-form is K\"ahler. \endproof

\hfill

%%%%%%%%%%%%%%%%%%%%%%%%%%%%%%%%%%%%%%%%%%%%%%%%%%%%%%%%%%%%
\proposition\label{_kahler_pote_on_Vaism_Proposition_}
Let $M$ be an LCK manifold with parallel Lee form $\theta$.
Assume that $\theta$ is exact, and let $r$ be the 
corresponding potential function. Then $r$ is the
K\"ahler potential for the K\"ahler form
$r\omega$.

\hfill

{\bf Proof:} Let $\Lie_{\theta^\sharp}$ be the operator of
Lie derivative along the vector field $\theta^\sharp$ dual
to $\theta$. Then $\Lie_{\theta^\sharp} \omega=0$
by \ref{_flow_along_theta_automo_Proposition_}.
Similarly,
\[ \Lie_{\theta^\sharp} r = dr \cntrct \theta^{\sharp} 
    = -r \theta \cntrct\theta^{\sharp} = -r.
\]
Therefore, $\Lie_{\theta^\sharp}(r \omega) = - r\omega$.
On the other hand, $r\omega$ is closed: $d(r\omega) =0$.
We obtain
\[  r\omega = - \Lie_{\theta^\sharp}(r \omega) = 
    d (r \omega\cntrct \theta^{\sharp}).
\]
Let $d^c= - I \circ d \circ I$ be the 
twisted de Rham differential. The function
$r$ is a K\"ahler potential for the form $r\omega$
if $r \omega = d d^c r$.
As we have seen above,
$r \omega = d(r \omega \cntrct \theta^{\sharp})$. Therefore, 
$r \omega = d d^c r$ is implied by 
\begin{equation} \label{_Kah_pote_1-fo_Equation_}
r \omega \cntrct \theta^{\sharp} = d^c r.
\end{equation}
To prove \eqref{_Kah_pote_1-fo_Equation_},
notice that $\omega \cntrct \theta^{\sharp} = I(\theta)$,
and $dr = r \theta$. Therefore, 
\begin{equation} \label{_d^c_r_Equation_}
   r \omega \cntrct \theta^{\sharp} = r I(\theta)
   = I(d r) = d^c r.
\end{equation}
This proves \eqref{_Kah_pote_1-fo_Equation_}.
\ref{_kahler_pote_on_Vaism_Proposition_}
is proven. \endproof

%%%%%%%%%%%%%%%%%%%%%%%%%%%%%%%%%%%%%%%%%%%%%%%%

\section{Einstein-Weyl LCK manifolds}

%%%%%%%%%%%%%%%%%%%%%%%%%%%%%%%%%%%%%%%%%%%%%%%%

In the Section we relate the definition and 
basic results on the geometry of Einstein-Weyl LCK manifolds.
We follow \cite{_Ornea:LCHK_}. 

\hfill

%%%%%%%%%%%%%%%%%%%%%%%%%%%%%%%%%%%%%%%%%%%%%%%%
\definition
Let $M$ be a Weyl manifold, 
$[g]$ its conformal class,
$\nabla$ its Weyl connection
and $R$ the Ricci curvature of $\nabla$,
$R = \Tr_{jl} \Theta_{ijk}^l$, where
$\Theta\in \Lambda^2(M)\otimes \Lambda^1(M) \otimes TM$
is the curvature of $\nabla$. Assume that the
symmetric part of $R$ is proportional to $[g]$.
Then $M$ is called {\bf an Einstein-Weyl manifold}. 

\hfill

The following important result was proven by P. Gauduchon.

\hfill

%%%%%%%%%%%%%%%%%%%%%%%%%%%%%%%%%%%%%%%%%%%%%%%%%%%%%%%%%%%%
\theorem\label{_E-W_theta_par_Theorem_}
\cite{_Gauduchon_19995_}
Let $M$ be a compact Einstein-Weyl manifold with
closed Lee form, and $g$ the Gauduchon metric on $M$
(\ref{_Gaud_metric_Definition_}). Denote by $\theta$
the corresponding Lee form. Assume that $\theta$
is not exact. Then 
\begin{description}
\item[(i)] $\theta$ is parallel with respect to the
Levi-Civita connection  associated with $g$.
\item[(ii)] The Ricci curvature of the Weyl connection
$\nabla$ vanishes.
\end{description}
\endproof

\hfill

%%%%%%%%%%%%%%%%%%%%%%%%%%%%%%%%%%%%%%%%%%%%%%%%
\remark 
By \ref{_E-W_theta_par_Theorem_} (i),
any Einstein-Weyl LCK manifold is also 
a Vaisman manifold.

\hfill

The following claim is quite obvious from the definitions. 

\hfill

%%%%%%%%%%%%%%%%%%%%%%%%%%%%%%%%%%%%%%%%%%%%%%%%
\claim\label{_Calabi_yau_covering_Einstein_}
Let $M$ be an LCK manifold equipped with a Gauduchon
metric, and $L$ the corresponding weight
bundle, equipped with a canonical
flat connection. Assume that $L$ has
trivial monodromy, that is, the form
$\theta$ is exact: $\theta = d(\zeta)$.
Let $e^{-\zeta}\omega$ be the corresponding
K\"ahler form associated with $\zeta$ as in
\ref{_theta_exact_confo_Riemann_Claim_}.
Then the following conditions are equivalent
\begin{description}
\item[(i)] The manifold $M$ is Einstein-Weyl
\item[(ii)] The K\"ahler metric $e^{-\zeta}\omega$ 
is Calabi-Yau, that is, the corresponding
Levi-Civita connection is Ricci-flat.
\end{description}

{\bf Proof:} Let $\nabla$ be the Weyl connection on $M$. Since $\nabla$
is torsion free and $\nabla$ preserves $e^{-\zeta}\omega$,
$\nabla$ is the Levi-Civita connection associated
with $e^{-\zeta}\omega$. The metric $e^{-\zeta}\omega$
is Calabi-Yau if and only if $\nabla$ is 
Ricci-flat. \endproof

\hfill

%%%%%%%%%%%%%%%%%%%%%%%%%%%%%%%%%%%%%%%%%%%%%%%%
\corollary
Let $M$ be a locally conformal hyperk\"ahler (LCHK)
manifold. Then $M$ is Einstein-Weyl. 

\hfill

{\bf Proof:} By \ref{_LCHK_w_trivial_monodro_Claim_},
the universal covering $\tilde M$ is hyperk\"ahler.
Therefore, $\tilde M$ is Calabi-Yau.
Now \ref{_Calabi_yau_covering_Einstein_}
implies that $M$ is Einstein-Weyl. 
\endproof

\hfill

Further on, we shall need the following
proposition.

\hfill

%%%%%%%%%%%%%%%%%%%%%%%%%%%%%%%%%%%%%%%%%%%%%%%%%%%%%%%%%%%%%%%%%%%%%%%%
\proposition\label{_can_cla_weight_on_Einstein_Proposition_}
Let $M$ be an Einstein-Weyl LCK manifold, and
$L$ its weight bundle. Then $L^n = K^{-1}$,
where $K$ is the canonical bundle of $M$
and $n = \dim_\C M$.

\hfill

{\bf Proof:} 
Let $\tilde M$ be a universal covering of $M$. By 
\ref{_Calabi_yau_covering_Einstein_}, $\tilde M$
is equipped with a Ricci-flat K\"ahler metric.
Therefore, the holonomy $\Hol(\tilde M)$
is contained in $SU(n)$. Since the monodromy
of $M$ preserves the conformal class
of the metric, we have $\Hol(M)\subset SU(n)\cdot \R^{>0}$,
where $\R^{>0}$ denotes the group of positive real numbers. 
By definition, the weight bundle corresponds to a 
representation $\det_\R(TM)^{\frac 1 {2n}}$, $n=\dim_\C M$ 
(see \ref{_weight_determi_Remark_}).
Therefore, the monodromy group of $L$
coincides with the quotient $\Hol(M)/\Hol(\tilde M)$
of full holonomy by the local holonomy.
Let $\alpha$ be a non-zero element of the monodromy 
group $G$ of $L$,
$G = \Hol(M)/\Hol(\tilde M)\subset \R^{>0}$.
Then $\alpha$
acts on $L$ as $(\alpha, l) \arrow \alpha l, \alpha\in \R^{>0}, l \in L$.

The action of $\Hol(M)$ on $K(M)$ factors through
\[ G = \Hol(M)/\Hol(\tilde M),\] because 
$\Hol(\tilde M)\subset SU(n).$
The quotient \[ \Hol(M)/\Hol(\tilde M)\subset \R^{>0}\] acts on 
$K(M) = \det(\Lambda^{1,0}(M))$ as
\[ (\alpha, \eta) \arrow \alpha^{-n} \eta,\ \ 
   \alpha\in \R^{>0}, \eta \in\det(\Lambda^{1,0}(M)).
\]
This relates the monodromy of $K(M)$ and the action
of holonomy.
The bundles $L^n$, $K^{-1}$ are flat, and their
monodromy is equal. Therefore, these bundles
are isomorphic.
\endproof

%%%%%%%%%%%%%%%%%%%%%%%%%%%%%%%%%%%%%%%%%%%%%%%%

\section{The form $\omega_0$ on Vaisman manifolds}
\label{_omega_0_Section_}

%%%%%%%%%%%%%%%%%%%%%%%%%%%%%%%%%%%%%%%%%%%%%%%%

In this Section we present some basic results and calculations on 
the geometry of Vaisman manifolds.

%%%%%%%%%%%%%%%%%%%%%%%%%%%%%%%%%%%%%%%%%%%%%%%%
\subsection{The form $\omega_0$: definition and eigenvalues}
%%%%%%%%%%%%%%%%%%%%%%%%%%%%%%%%%%%%%%%%%%%%%%%%

Let $M$ be an LCK manfold. 
Consider the form 
\[ 
  \omega_0:= d^c \theta
\] 
on $M$. We have $d^c = \frac{\6 - \bar\6}{\1}$
Therefore,  
\[ \omega_0 = \frac{\6 - \bar\6}{\1}\theta.
\]
Write the Hodge decomposition of $\theta$ as 
$\theta = \theta^{1,0} + \theta^{0,1}$.
Since $\theta$ is closed, we have 
\[ \6 \theta^{1,0} = \bar\6 \theta^{0,1}=0, 
    \ \ \6 \theta^{0,1} = - \bar\6\theta^{1,0}.
\]
This implies
\begin{equation}\label{_omega_0_6_Equation_} 
   \omega_0  = -2\1 \6\theta^{0,1}
\end{equation}

%%%%%%%%%%%%%%%%%%%%%%%%%%%%%%%%%%%%%%%%%%%%%%%%%%%%%%%%%%%%
\proposition \label{_omega_0_eigenvalues_Proposition_}
Let $M$ be a Vaisman manifold, that is, 
an LCK manfold with parallel Lee form $\theta$.
Consider the form $\omega_0:= d^c \theta$.
Chose an orthonormal basis 
\[ \xi_1, ..., \xi_{n-1}, \sqrt 2 \cdot\theta^{1,0} \]
in $T^{1,0}M$,  where $\theta^{1,0}$ is the 
$(0,1)$-part of $\theta$,\footnote{Since $|\theta|=1$, we have 
$|\theta^{1,0}|=\frac{1}{\sqrt 2}$.}
and let 
\[ \omega= \1(\xi_1 \wedge \bar\xi_1 + 
   \xi_2 \wedge \bar\xi_2 + ... + 
   \xi_{n-1} \wedge \bar\xi_{n-1} + 2 \theta^{1,0}\wedge \theta^{0,1})
\]
be the Hermitian form on $M$.
Then 
\begin{equation}\label{_omega_0_coord_Equation_}
\omega_0 = \1(\xi_1 \wedge \bar\xi_1 + 
   \xi_2 \wedge \bar\xi_2 + ... + 
   \xi_{n-1} \wedge \bar\xi_{n-1}).
\end{equation}
In particular, all eigenvalues of $\omega_0$ are 
positive except the one corresponding to $\theta$, 
which is equal zero.

\hfill

{\bf Proof:} Passing to a covering, we can always assume
that $\theta$ is exact. Let $r$ be the potential of $M$.
By \eqref{_d^c_r_Equation_}, we have
\begin{equation} \label{_d^c_r_theta_Kahler_Equation_}
d^c(r\theta) = r \omega.
\end{equation}
On the other hand,
\begin{equation} \label{_d^c_r_theta_preci_Equation_}
d^c(r\theta) = r d^c \theta + d^c r \wedge \theta 
             = r d^c \theta + r I(\theta) \wedge \theta.
\end{equation}
Comparing \eqref{_d^c_r_theta_Kahler_Equation_} 
and \eqref{_d^c_r_theta_preci_Equation_}, we find
$\omega = d^c\theta + \theta \wedge I(\theta)$,
and 
\begin{equation}\label{_omega_0_via_theta_wedge_theta_c_Equation_}
 \omega_0 = \omega - 2\1\theta^{1,0}\wedge \theta^{0,1}.
\end{equation}
This proves 
\ref{_omega_0_eigenvalues_Proposition_}.
\endproof

\hfill

The following claim is obvious.

\hfill

%%%%%%%%%%%%%%%%%%%%%%%%%%%%%%%%%%%%%%%%%%%%%%%%
\claim\label{_omega_0_exact_Claim_}
The form $\omega_0$ is exact. 

\hfill

{\bf Proof:} By definition, we have
\begin{equation}\label{_omega_0_defi_Equation_}
\omega_0 = - I d(I(\theta)). 
\end{equation}
On the other hand,
$\omega_0$ is a (1,1)-form, hence
$I(\omega_0) = \omega_0$. Comparing this with
\eqref{_omega_0_defi_Equation_}, we find
\[ 
\omega_0 = - d(I(\theta)).
\]
\endproof

%%%%%%%%%%%%%%%%%%%%%%%%%%%%%%%%%%%%%%%%%%%%%%%%%%%%%%%%%%%%
\subsection{The form $\omega_0$ and the canonical foliation}
\label{_cano_foli_Subsection_}
%%%%%%%%%%%%%%%%%%%%%%%%%%%%%%%%%%%%%%%%%%%%%%%%%%%%%%%%%%%%

Let $M$ be a Vaisman manifold, that is, 
an LCK manifold with a parallel Lee form $\theta$. 
Consider a 2-dimensional real foliation $\Xi \subset TM$
generated by $\theta^\sharp, I(\theta^\sharp)$. 
The vector field $\theta^\sharp$ is holomorphic
(\ref{_flow_along_theta_automo_Proposition_}). 
Therefore, $\Xi$ is integrable and holomorphic.

\hfill

%%%%%%%%%%%%%%%%%%%%%%%%%%%%%%%%%%%%
\definition
The foliation $\Xi$ is called {\bf the canonical
foliation} of a Vaisman manifold $M$. 

\hfill

The canonical foliation is related to the form
$\omega_0$ in the following way. 

\hfill

%%%%%%%%%%%%%%%%%%%%%%%%%%%%%%%%%%%%%%%%%%%%%%%%%%%%%%%%%%%%%%%%%%%%%%%%
\proposition\label{_Kahler_form_on_space_of_leaves_Proposition_}
Let $M$ be a Vaisman manifold, $\Xi$ the 
foliation defined above and $\omega_0 = d^c \theta$ the
standard $(1,1)$-form. Assume that the space of leaves
of $\Xi$ is well defined, and let $f:\; M \arrow Y$
the the corresponding quotient map. 
Then $Y$ is equipped with a natural K\"ahler
form $\omega_Y$, in such a way that $f^*\omega_Y = \omega_0$.

\hfill

{\bf Proof:}
To check that $\omega_0 = f^* \omega_Y$ for some 2-form $\omega_Y$,
we need to show that 
\begin{equation}\label{_theta_sharp_into_omega_0_Equation_}
\omega_0 \cntrct \theta^\sharp = \omega_0 \cntrct I(\theta^\sharp) =0
\end{equation}
and
\begin{equation}\label{_theta_sharp_Lie_deri_omega_0_Equation_}
\Lie_{\theta^\sharp}\omega_0 = \Lie_{I(\theta^\sharp)}\omega_0=0.
\end{equation}
The equation \eqref{_theta_sharp_into_omega_0_Equation_}
follows  immediately from \ref{_omega_0_eigenvalues_Proposition_}.
Using the Cartan formula for the the Lie derivative and
\eqref{_theta_sharp_into_omega_0_Equation_}, we obtain
\[
  \Lie_{\theta^\sharp}\omega_0 = (d\omega_0)\cntrct \theta^\sharp =0
\]
and 
\[
  \Lie_{I(\theta^\sharp)}\omega_0 = (d\omega_0)\cntrct I(\theta^\sharp) =0.
\]
This implies that $\omega_0 = f^* \omega_Y$ for some 2-form
$\omega_Y$ on $Y$. Since $\omega_0$ has $\dim_\C Y=n-1$
positive eigenvalues, the form $\omega_Y$ is positive
definite. Therefore, $\omega_Y$ is a K\"ahler form on $Y$.
\endproof

\hfill

The form $\omega_0$ is exact (\ref{_omega_0_exact_Claim_}).
The following proposition immediately follows from
this observation.

\hfill

%%%%%%%%%%%%%%%%%%%%%%%%%%%%%%%%%%%%%%%%%%%%%%%%
\proposition\label{_subva_Vaisman_Proposition_}
Let $M$ be a compact Vaisman manifold,
and $X\subset M$ a closed complex subvariety.
Then $X$ is tangent to the canonical foliation $\Xi$
in all its smooth points. In particular, if $X$ is smooth,
then $X$ is also a Vaisman manifold.

\hfill

{\bf Proof:} The form $\omega_0$ is {\bf positive}
in the sense of distribution theory;
that is, $\omega_0$ is a real $(1,1)$-form
with non-negative eigenvalues. It is well known that
\begin{equation}\label{_int_of_omega_0_Equation_}
\int_X \omega_0^k\geq 0
\end{equation}
for all complex subvarieties
$X\subset M$, $\dim_\C X =k$, and all positive forms $\omega_0$. 
Moreover, the integral \eqref{_int_of_omega_0_Equation_}
vanishes only if $X$ is tangent to the
null-space foliation of $\omega_0$. 

Since $\omega_0$ is exact, the integral 
\eqref{_int_of_omega_0_Equation_}
vanishes. Therefore, $X$ is tangent to the
null-space foliation of $\omega_0$. As we have seen above,
the null-space of $\omega_0$ is $\Xi$. This implies
$X$ is tangent to the canonical foliation.

If $X$ is smooth, it is clearly an LCK manifold.
To prove that $X$ is a Vaisman submanifold, we use
the following theorem of Kamishima and Ornea
(\cite{_Kamishima_Ornea_}). 

\hfill

%%%%%%%%%%%%%%%%%%%%%%%%%%%%%%%%%%%%%%%%%%%%%%%%
\theorem\label{_Kamishima_Ornea_Theorem_}
\cite{_Kamishima_Ornea_}
Let $X$ be a 
compact LCK manifold admitting a conformal holomorphic
flow which is not conformal equivalent to isometry.
Then $X$ is a Vaisman manifold.

\endproof

\hfill

Consider the holomorphic flow $\psi_t$ associated with the Lee
field on $M$ (\ref{_flow_along_theta_automo_Proposition_}).
Since $X$ is tangent to $\Xi$, $\psi_t$ preserves $X$.
It is trivial to check that $\psi_t$ is not 
conformal equivalent to isometry.
Applying \ref{_Kamishima_Ornea_Theorem_},
we find that $X$ is Vaisman. 
\ref{_subva_Vaisman_Proposition_} is proven. \endproof

\hfill

For smooth manifolds, \ref{_subva_Vaisman_Proposition_}
is proven by K. Tsukada, see \cite{_Tsukada:subva_}.

%%%%%%%%%%%%%%%%%%%%%%%%%%%%%%%%%%%%%%%%%%%%%%%%
\subsection{Curvature of a weight bundle}
\label{_weight_curv_Subsection_}
%%%%%%%%%%%%%%%%%%%%%%%%%%%%%%%%%%%%%%%%%%%%%%%%

Let $M$ be an LCK manifold and $L$ its weight bundle.
By construction, $L$ is a complex vector bundle
 equipped with a flat connection
${}^L\nabla$ and a nowhere degenerate section $\lambda$,
such that
\[  {}^L\nabla(\lambda) = \lambda\otimes\left(-\frac{1}{2}\theta\right).
\]
A $(0,1)$-part of ${}^L\nabla$ gives a holomorphic structure
on $L$. Throughout this Subsection, we shall
consider $L$ as a holomorphic bundle. 
Consider a Hermitian structure $g_L$ on $L$,
defined in such a way that $|\lambda|_{g_L}=1$.

\hfill

%%%%%%%%%%%%%%%%%%%%%%%%%%%%%%%%%%%%%%%%%%%%%%%%
\theorem\label{_curva_L_omega_0_Theorem_}
Let $M$ be an LCK manifold, and $L$ the 
corresponding weight bundle equipped with a 
holomorphic and a Hermitian structure as above.
Let ${}^C\nabla$ be the standard Hermitian connection on $L$
(so-called Chern connection), and $C$ its
curvature. Then $C= -2 \1\omega_0$, where 
$\omega_0 = d^c \theta$ is the 
standard 2-form on $M$.

\hfill

{\bf Proof:} By definition of the Chern connection, we have
\begin{equation} \label{_Chern_conne_on_lambda_Equation_}
\begin{aligned}
{}^C\nabla^{0,1}(\lambda) = {}^L\nabla^{0,1}(\lambda) &=
\lambda \otimes \left(-\frac{1}{2}\theta^{0,1}\right), \text{\ \ and\ \ } \\
( {}^C\nabla^{1,0}(\lambda), \lambda)_H &= - (\lambda, {}^C\nabla^{0,1}(\lambda))_H,
\end{aligned}
\end{equation}
where $(\cdot,\cdot)_H$ is the Hermitian form on $L$.
{}From \eqref{_Chern_conne_on_lambda_Equation_}, we obtain
\begin{equation}\label{_Chern_conne_on_weight_Equation_}
{}^C\nabla= \nabla_{tr} + \frac{\theta^{0,1}-\theta^{1,0}}{2},
\end{equation}
where $\nabla_{tr}$ is the trivial connection on $L$
fixing $\lambda$. From \eqref{_Chern_conne_on_weight_Equation_},
we obtain
\begin{equation}\label{_weigh_conne_curv_Equation_}
C = \frac{\6\theta^{0,1} - \bar\6\theta^{0,1}}{2}.
\end{equation}
On the other hand,  $\theta$ is closed, and therefore
$\6\theta^{0,1} = - \bar\6\theta^{0,1}$.
Comparing this with \eqref{_weigh_conne_curv_Equation_},
we obtain 
\[ C = \6 \theta^{0,1} = -2 \1\omega_0.\]
(the last equation holds by \eqref{_omega_0_6_Equation_}). 
We proved \ref{_curva_L_omega_0_Theorem_}. \endproof

%%%%%%%%%%%%%%%%%%%%%%%%%%%%%%%%%%%%%%%%%%%%%%%%%%%%%%%%%%%%

\section{K\"ahler geometry of the form $\omega_0$}

%%%%%%%%%%%%%%%%%%%%%%%%%%%%%%%%%%%%%%%%%%%%%%%%%%%%%%%%%%%%

The form $\omega_0$ behaves, in many ways, as a surrogate K\"ahler form
on $M$. In this Section we prove the $\omega_0$-version of the
Kodaira relations and apply it to obtain the standard
identities for Laplace operators with coefficients in a bundle.
This leads to Kodaira-Nakano-type vanishing theorems.

%%%%%%%%%%%%%%%%%%%%%%%%%%%%%%%%%%%%%%%%%%%%%%%%%%%%%%%%%%%%
\subsection{The $SL(2)$-triple $L_0, \Lambda_0, H_0$}
%%%%%%%%%%%%%%%%%%%%%%%%%%%%%%%%%%%%%%%%%%%%%%%%%%%%%%%%%%%%

In this Subsection, we study
the Lefschetz-type $SL(2)$-action associated with $\omega_0$.
Let $M$ be a Vaisman manifold and $\omega_0$ the standard
2-form (Section \ref{_omega_0_Section_}). Denote by
$L_0$ the operator $\eta \arrow \eta \wedge \omega_0$,
and by $\Lambda_0$ the Hermitian adjoint operator.
Using the coordinate expression
of $\omega_0$ given in \ref{_omega_0_eigenvalues_Proposition_},
we find that $L_0$, $\Lambda_0$, $H_0:= [ L_0, \Lambda_0]$
form an $SL(2)$-triple (exactly the same argument
is used in the proof of Lefschetz theorem via
the $SL(2)$-action on cohomology, see 
\cite{_Griffi_Harri_}).

Locally, the operator $H_0$ can be expressed 
as follows. Let 
\[ \xi_1, \xi_2, ..., \xi_{n-1}, 
   \sqrt 2 \theta^{1,0}, 
   \bar\xi_1, \bar\xi_2, ..., 
   \bar\xi_{n-1}, \sqrt 2 \theta^{0,1}\in \Lambda^1(M)
\]
be an orthonormal frame in the bundle of forms
(\ref{_omega_0_eigenvalues_Proposition_}). 
Consider a monomial
\begin{equation}\label{_monomial_H_0_Equation_}
   \lambda = \xi_{i_1} \wedge \xi_{i_2}\wedge ... \wedge \xi_{i_k}
   \wedge \bar\xi_{i_{k+1}} \wedge \bar\xi_{i_{k+2}}\wedge ... 
   \wedge \xi_{p} \wedge R
\end{equation}
where $R$ is a monomial of $\theta^{0,1}, \theta^{1,0}$.
Then 
\begin{equation}\label{_H_0_explicitly_Equation_}
H_0(\lambda) = (p - n +1) \lambda.
\end{equation}
where  $p$ is the number of $\xi_{i_l}, \bar\xi_{i_k}$
in $\lambda$ and $n = \dim_\C M$.

The equation \eqref{_H_0_explicitly_Equation_}
is proved in the same way as the explicit form of the
operator $H$ in Hodge theory (\cite{_Griffi_Harri_}).

%%%%%%%%%%%%%%%%%%%%%%%%%%%%%%%%%%%%%%%%%%%%%%%%%%%%%%%%%%%%%%%%%%%%%%%%
\subsection{Kodaira identities for the differential $\6_0$}
\label{_Koda_for_6_0_Subsection_}
%%%%%%%%%%%%%%%%%%%%%%%%%%%%%%%%%%%%%%%%%%%%%%%%%%%%%%%%%%%%%%%%%%%%%%%%

Let $M$ be a Vaisman manifold (an LCK manifold with
a parallel Lee form). Consider the $\omega_0$-Lefschetz
triple $L_0$, $\Lambda_0$, $H_0$ acting on $\Lambda^*(M)$
as above, and let
\[ \Lambda^*(M) = \bigoplus_i \Lambda^*_i(M)
\]
be the weight decomposition associated with $H_0$,
in such a way that the monomial 
\eqref{_monomial_H_0_Equation_} has weight $p$.
Clearly, 
\begin{equation}\label{_H_0-decompo_weights_Equation_}
\Lambda^p(M) = \Lambda^p_{p-2}(M) \oplus \Lambda^p_{p-1}(M)
 \oplus \Lambda^p_{p}(M).
\end{equation}
Denote by $d_0$ the weight 1 component of 
the de Rham differential:
\begin{equation}\label{_diffe_d_0_weight_Equation_}
   d_0:\; \Lambda^p_{i}(M)\arrow \Lambda^{p+1}_{i+1}(M).
\end{equation}

\hfill

%%%%%%%%%%%%%%%%%%%%%%%%%%%%%%%%%%%%%%%%%%%%%%%%
\remark\label{_d_0_from_omega_0_etc_Remark_}
By \eqref{_omega_0_via_theta_wedge_theta_c_Equation_}, we have
$d(\omega_0) = \omega_0 \wedge \theta$. Therefore,
$d_0(\omega_0)=0$. Similarly, one can check 
by elementary computations that
$d_0 \theta=0$ and $d_0(I(\theta))=0$.

\hfill

%%%%%%%%%%%%%%%%%%%%%%%%%%%%%%%%%%%%
\remark
The differential $d_0$ satisfies the Leibniz rule.

\hfill

The $(1,0)$-part and $(0,1)$-part
of $d_0$ are denoted by $\6_0$, $\bar\6_0$ as
usual. 

\hfill

%%%%%%%%%%%%%%%%%%%%%%%%%%%%%%%%%%%%%%%%%%%%%%%%
\definition
The operator $d_0$ is called {\bf the 
$\omega_0$-de Rham differential}, the operators $\6_0$, $\bar\6_0$ 
{\bf the $\omega_0$-Dolbeault operators}.

\hfill

%%%%%%%%%%%%%%%%%%%%%%%%%%%%%%%%%%%%%%%%%%%%%%%%
\theorem \label{_Kodaira_for_d_0_Theorem_}
(Kodaira identities for $\6_0$, $\bar\6_0$).
Let $M$ be a Vaisman manifold,
and $\6_0$, $\bar\6_0$, $L_0$, $\Lambda_0$ the
operators defined above. Then
\begin{equation}\label{_Kodaira_relations_6_0_Equation_}
\begin{aligned}
\ & [\Lambda_0, \6_0] = \1 \bar\6_0^*,  \ \ \ 
 [L_0, \bar\6_0] = - \1 \6_0^*,  \\
\ &  [\Lambda_0, \bar\6_0^*] = - \1 \6_0, \ \ \ 
 [L_0, \6_0^*] = \1 \bar\6_0, 
\end{aligned}
\end{equation}
where $\6_0^*$, $\bar\6_0^*$ are the Hermitian adjoint 
operators to $\6_0$, $\bar\6_0$.

\hfill

{\bf Proof:} 
Further on, we shall prove
\begin{equation}\label{_Kodaira_d_d^c_Equation_}
[L_0, d_0^*] = - I \circ d_0 \circ I.
\end{equation}
Taking the $(0,1)$ and $(1,0)$-parts of 
\eqref{_Kodaira_d_d^c_Equation_}, we obtain the 
bottom line of \eqref{_Kodaira_relations_6_0_Equation_}.
Taking Hermitian adjoint, we obtain the top line of
\eqref{_Kodaira_relations_6_0_Equation_}. Therefore,
\eqref{_Kodaira_d_d^c_Equation_} implies 
the Kodaira relations for $\6_0$.

We prove \eqref{_Kodaira_d_d^c_Equation_} 
in an algebraic fashion similar to the proof
of Kodaira relations in HKT geometry (\cite{_Verbitsky:HKT_}).
Consider $L_0$, $d_0$, etc. as operators on the algebra
of differential forms. A. Grothendieck 
gave a general recursive definition of
differential operators on an algebra
(see, e.g., \cite{_Verbitsky:HKT_}).\footnote{An 
$n$-th order differential operator
on an algebra $\Lambda^*(M)$ can have any order 
{}from 0 to $n$ as a differential operator in the usual sense.
To avoid confusion, we shall write ``algebraic differential
operator'' speaking of differential operators, on the algebra $\Lambda^*(M)$,
in the sense of Grothendieck.}
Then $L_0$ is a 0-th order algebraic differential operator,
and $d_0^*$ a second order algebraic differential operator
on $\Lambda^*(M)$
(this is straightforward; see the full
argument in \cite{_Verbitsky:HKT_}). Therefore, the commutator
$[L_0, d_0^*]$ is a first order algebraic differential operator.
Since $- I \circ d_0 \circ I$ satisfies the Leibniz rule,
it is also a first order algebraic differential operator. 
A first order algebraic differential operator $D$ satisfies
\begin{equation}\label{_diffe_ope_Equation_} 
 D(ax) = (-1)^{\deg D \deg a} a D(x) + D(a) x - D(1) ax. 
\end{equation}
This is clear from the definition; see,
again, \cite{_Verbitsky:HKT_}.
{}From \eqref{_diffe_ope_Equation_}, 
we obtain that the first order algebraic differential operator
is determined by the values taken on any
set of generators of the algebra.
Since both sides of \eqref{_Kodaira_d_d^c_Equation_}
are  first order algebraic differential operators, it suffices 
to check that they are equal on some subspace 
$V\subset \Lambda^*(M)$ which generates 
the algebra $\Lambda^*(M)$. 

On 0-forms, \eqref{_Kodaira_d_d^c_Equation_} is clear:
\[ d_0^*(f \omega_0) = \omega_0 \cntrct (d_0 f)^\sharp = - I(d_0 f) 
\]
where ${v}^\sharp$, as usually, 
denotes the vector field associated
with a form $v$. 

We are going to construct a subspace $V \subset \Lambda^1(M)$
generating $\Lambda^1(M)$ over $C^\infty(M)$ such that
the operators on both sides of equation 
\eqref{_Kodaira_d_d^c_Equation_} are equal on $V$.
As we have mentioned above, this is sufficient
for the proof of \ref{_Kodaira_for_d_0_Theorem_}.

The statement of \ref{_Kodaira_for_d_0_Theorem_} is local.
Passing to an open neighbourhood
if necessary, we may assume that the space $Y$
of leaves of $\Xi$ is well defined. Let $f:\; M \arrow Y$
be the corresponding quotient map. By 
\ref{_Kahler_form_on_space_of_leaves_Proposition_},
the manifold $Y$ is equipped with a natural K\"ahler
form $\omega_Y$, in such a way that $f^* \omega_Y = \omega_0$.
Let $d_Y$ be the de Rham operator on $Y$,
and $d_Y^*$ its Hermitian adjoint. 
Denote by $L_Y$ the Hodge operator on $Y$,
$L_Y(\eta) = \eta \wedge\omega_Y$. Since $Y$
is K\"ahler, the usual Kodaira identity holds:
\begin{equation}\label{_Kodaira_d_d^c_on_Y_Equation_}
[L_Y, d_Y^*] = - I \circ d_Y \circ I,
\end{equation}
Lifting \eqref{_Kodaira_d_d^c_on_Y_Equation_}
to $M$, we obtain that \eqref{_Kodaira_d_d^c_Equation_}
holds on all forms $\eta = f^* \eta_Y$ which 
are obtained as a pullback. 

Let us check \eqref{_Kodaira_d_d^c_Equation_}
on the 2-dimensional space 
$\langle \theta, I(\theta)\rangle\subset \Lambda^1(M)$
generated by $\theta$ and $I(\theta)$. 
We have $- I \circ d_0\circ I (I \theta)=0$ because $\theta$ is 
$d_0$-closed (\ref{_d_0_from_omega_0_etc_Remark_}).
The same argument proves
\[ 
d_0^*(\omega^k \wedge I(\theta))= * d_0(\omega^{n-k} \wedge \theta)=0
\]
because $d_0(\theta) = d_0 (\omega)=0$ 
(see \ref{_d_0_from_omega_0_etc_Remark_}).
Therefore, the operators 
$[L_0, d_0^*]$ and $- I \circ d_0\circ I $
vanish on $I(\theta)$.

Similarly, 
we have $- I \circ d_0\circ I (\theta) = \omega_0$, and
$[L_0, d_0^*] (\theta)=0$.

We have shown that the operators 
$[L_0, d_0^*]$ and $- I \circ d_0\circ I $
are equal on the space $V \subset \Lambda^1(M)$
generated by $\theta, I(\theta)$ and the pullbacks of
differential forms from $Y$. Clearly, 
$C^\infty (M) \cdot V = \Lambda^1 (M)$.
Therefore, the first order algebraic differential operators
$[L_0, d_0^*]$ and $- I \circ d_0\circ I $ are equal
on a set of generators of $\Lambda^*(M)$.
This proves \eqref{_Kodaira_d_d^c_Equation_}.
\ref{_Kodaira_for_d_0_Theorem_} is proven. 
\endproof

\hfill

A similar argument proves the following

\hfill

%%%%%%%%%%%%%%%%%%%%%%%%%%%%%%%%%%%%%%%%%%%%%%%%
\claim\label{_d_0^2=0_Claim_}
Let $M$ be a Vaisman manifold, and $d_0$, $\6_0$, $\bar\6_0$ the
differential operators defined above. Then 
\[ d_0^2 = \6_0^2 = \bar\6_0^2=0.
\]
{\bf Proof:}
Since $\6_0^2$ and $\bar\6_0^2$ are $(2,0)$ and $(0,2)$-parts
of $d_0^2$, it suffices to show that $d_0^2=0$.
We use the same algebraic argument as in the proof of
\ref{_Kodaira_for_d_0_Theorem_}. The anticommutator
$d_0^2= \{d_0, d_0\}$ is a first order algebraic differential
operator, because it is a supercommutator of two 
first order algebraic differential operators.
Therefore, it suffices to prove that
$d_0^2=0$ on some set of generators of $\Lambda^*(M)$.

All 1-forms have weight $\leq 1$ with respect to
the decomposition \eqref{_H_0-decompo_weights_Equation_}.
Therefore, the restriction
$d\restrict{\Lambda^0(M)}$ has components
of weight 0 and 1 only (no weight 2 component).
Therefore, the weight 2 component of 
$d^2\restrict{\Lambda^0(M)}$ is equal to 
$d_0^2$. Since $d^2=0$, we have 
\[ d_0^2\restrict{\Lambda^0(M)}=0. 
\]
Passing to a local neighbourhood, we may always assume that
the space $Y$ of leaves of $\Xi$ is well defined.
Clearly, $d_0=d$ on all forms lifted from $Y$.
Therefore, $d_0^2=d^2=0$ on such forms.

Finally, $d_0=0$ on the space generated by 
$\theta, I(\theta)$. Hence $d_0^2=0$ on this space.
We obtain that $d_0^2=0$ on a set of generators
of the algebra $\Lambda^*(M)$. This implies 
$d_0^2=0$. We proved \ref{_d_0^2=0_Claim_}.
\endproof

%%%%%%%%%%%%%%%%%%%%%%%%%%%%%%%%%%%%%%%%%%%%%%%%%%%%%%%%%%%%
\subsection{Laplace operators with coefficients in a bundle}
\label{_Laplace_in_bundle_Subsection_}
%%%%%%%%%%%%%%%%%%%%%%%%%%%%%%%%%%%%%%%%%%%%%%%%%%%%%%%%%%%%

Let $M$ be a Vaisman manifold, $B$ a holomorphic vector bundle
equipped with a Hermitian form, and 
\[ 
   \nabla:\; \Lambda^p(M)\otimes B \arrow  \Lambda^{p+1}(M)\otimes B
\]
the Chern connection
on $B$. Denote by $d_0$ the weight 1 component of $\nabla$,
taken with respect to the 
$H_0$-action (see \eqref{_diffe_d_0_weight_Equation_}). 
Let $\6_0$, $\bar\6_0$ be the $(1,0)$- and $(0,1)$-parts
of $d_0$, and $\6_0^*$, $\bar\6_0^*$ the Hermitian
adjoint operators. 

\hfill

%%%%%%%%%%%%%%%%%%%%%%%%%%%%%%%%%%%%%%%%%%%%%%%%
\proposition \label{_Kodaira_for_d_0_coeff_Proposition_}
(Kodaira identities for $\6_0$, $\bar\6_0$ with coefficients
in a Hermitian bundle).
Let $M$ be a Vaisman manifold, $B$ a Hermitian holomorphic bundle
and 
\[ \6_0, \bar\6_0, L_0, \Lambda_0:\; 
   \Lambda^*(M)\otimes B  \arrow \Lambda^*(M)\otimes B 
\]
the operators defined above. Then
\begin{equation}\label{_Kodaira_relations_6_0_w_bun_Equation_}
\begin{aligned}
& [\Lambda_0, \6_0] = \1 \bar\6_0^*,  \ \ \ 
[L_0, \bar\6_0] = - \1 \6_0^*,  \\
 & [\Lambda_0, \bar\6_0^*] = - \1 \6_0, \ \ \ 
[L_0, \6_0^*] = \1 \bar\6_0.
\end{aligned}
\end{equation}

{\bf Proof:} The proof of 
\ref{_Kodaira_for_d_0_coeff_Proposition_} is essentially
the same as the proof of the Kodaira identities
with coefficients in a bundle on a K\"ahler manifold. We deduce
\ref{_Kodaira_for_d_0_coeff_Proposition_} from the usual
(coefficient-less) Kodaira identities,
proven in \ref{_Kodaira_for_d_0_Theorem_}.

Write the Chern connection in $B$ as
\begin{equation}\label{_Chern_forms_Equation_}
\nabla = \6^{tr} + \bar\6^{tr}+\eta+\bar\eta,
\end{equation}
where $\6^{tr} + \bar\6^{tr}$ is the trivial
connection fixing the Hermitian structure, and 
$\eta\in \Lambda^{1,0}(M)\otimes \End B$ a (1,0)-form.
Let $\eta_0:= \Pi(\eta)$ be the projection of
$\eta$ to $\Lambda^1_0(M)\otimes \End B\subset 
\Lambda^1(M)\otimes \End B$. We consider $\eta$, 
$\eta_0$ as operators on differential forms.
Denote by $\eta^*$, $\eta_0^*$ the Hermitian
adjoint operators. Using the coordinate expression
of $\omega_0$ given in \ref{_omega_0_eigenvalues_Proposition_}, 
we find
\begin{equation}\label{_Kodaira_for_eta_Equation_}
\begin{aligned}
\ & [\Lambda_0, \eta_0] = \1 \bar\eta_0^*,  \ \ \ 
[L_0, \bar\eta_0] = - \1 \eta_0^*,  \\
 \ & [\Lambda_0, \bar\eta^*_0] = - \1 \eta_0, \ \ \ 
[L_0, \eta^*_0] = \1 \bar\eta_0.
\end{aligned}
\end{equation}
Adding  \eqref{_Kodaira_for_eta_Equation_}
and \eqref{_Kodaira_relations_6_0_Equation_}
(which holds for the trivial connection
$\nabla^{tr}=\bar\6^{tr} + \bar\6^{tr}$)
termwise, we obtain 
\eqref{_Kodaira_relations_6_0_w_bun_Equation_}.
This proves 
\ref{_Kodaira_for_d_0_coeff_Proposition_}.
\endproof

\hfill

%%%%%%%%%%%%%%%%%%%%%%%%%%%%%%%%%%%%%%%%%%%%%%%%
\theorem\label{_Lapla_compa_Theorem_}
Let $M$ be a Vaisman manifold, $B$ a Hermitian holomorphic bundle
and 
\[ \6_0, \bar\6_0, L_0, \Lambda_0:\; 
   \Lambda^*(M)\otimes B  \arrow \Lambda^*(M)\otimes B 
\]
the operators defined above.
Consider the Laplacians
$\Delta_{\6_0}= \6_0\6_0^* +\6_0^*\6_0$,
$\Delta_{\bar\6_0}= \bar\6_0\bar\6_0^* +\bar\6_0^*\bar\6_0$.
Then
\begin{equation}\label{_Lapla_compa_Equation_}
\Delta_{\6_0}-\Delta_{\bar\6_0} = \1 [ \Theta_B, \Lambda_0],
\end{equation}
where $\Theta_B:\; \Lambda^p(M) \arrow \Lambda^{p+2}(M)\otimes B$
is the curvature operator of $B$.

\hfill

{\bf Proof:} \ref{_Lapla_compa_Theorem_} is a formal
consequence of the Kodaira relations
\eqref{_Kodaira_relations_6_0_w_bun_Equation_}
(see e.g. \cite{_Griffi_Harri_} for a detailed proof).
\endproof

%%%%%%%%%%%%%%%%%%%%%%%%%%%%%%%%%%%%%%%%%%%%%%%%%%%%%%%%%%%%
\subsection{Serre's duality for $\bar\6_0$-cohomology}
%%%%%%%%%%%%%%%%%%%%%%%%%%%%%%%%%%%%%%%%%%%%%%%%%%%%%%%%%%%%

Further on, we shall use the following version of Serre's
duality. Since the operator $\bar\6_0$ is not elliptic,
the $\bar\6_0$-cohomology can be infinite-dimensional.
Therefore, it is more convenient to state the Serre's duality
as an isomorphism of vector spaces. 

\hfill

%%%%%%%%%%%%%%%%%%%%%%%%%%%%%%%%%%%%%%%%%%%%%%%%%%%%%%%%%%%%
\theorem\label{_Serre_s_for_6_0_Theorem_}
(Serre's duality for $\bar\6_0$-cohomology)
Let $M$ be a compact Vaisman manifold, 
$\dim_\C M =n$, $B$ a holomorphic Hermitian bundle,
and $K$ the canonical bundle. Then there exists an isomorphism
of the spaces of $\bar\6_0$-harmonic forms
with coefficients in $B$, $B^*\otimes K$
\[
{\cal H}_{\bar \6_0}^p(B) \cong 
\overline{{\cal H}_{\bar \6_0}^{n-p}(B^*\otimes K)},
\]
where $\overline V$ denotes the complex conjugate
vector space to $V$, that is, the same real space with the 
opposite complex structure.

\hfill

{\bf Proof:} Consider the Hodge $*$ operator acting on 
the differential forms with coefficients in a bundle.
Given a ${\bar \6_0}$-harmonic form $\eta\in \Lambda^{0,p}(M, B)$,
the form $*\eta$ is a ${\6_0}$-harmonic $(n-p,n)$-form
with coefficients in $\bar B^*$. This is clearly the same as
${\6_0}$-harmonic $(n-p, 0)$-form with coefficients in
$\bar K \otimes \bar B^*$. Taking a complex conjugate,
we obtain that $\overline{*\eta}$ can be considered
as a ${\bar \6_0}$-harmonic form with coefficients
in $B^*\otimes K$. This proves \ref{_Serre_s_for_6_0_Theorem_}.
\endproof

%%%%%%%%%%%%%%%%%%%%%%%%%%%%%%%%%%%%%%%%%%%%%%%%%%%%%%%%%%%%

\section{Vanishing theorems for Vaisman manifolds}

%%%%%%%%%%%%%%%%%%%%%%%%%%%%%%%%%%%%%%%%%%%%%%%%%%%%%%%%%%%%

%%%%%%%%%%%%%%%%%%%%%%%%%%%%%%%%%%%%%%%%%%%%%%%%%%%%%%%%%%%%
\subsection{Vanishing theorem for the differential $\6_0$}
\label{_vani_for_6_0_Subsection_}
%%%%%%%%%%%%%%%%%%%%%%%%%%%%%%%%%%%%%%%%%%%%%%%%%%%%%%%%%%%%

As one does in the proof of Kodaira-Nakano-type vanishing theorem, the
Kodaira relations \eqref{_Lapla_compa_Equation_} can be used to obtain
various vanishing results for $\bar\6_0$-cohomology
of holomorphic vector bundles. We do not need the whole
spectrum of vanishing theorems in this paper;
they can be stated and proven in the same way as the
usual vanishing theorems. For our present purposes,
we need only the following result.

\hfill

%%%%%%%%%%%%%%%%%%%%%%%%%%%%%%%%%%%%%%%%%%%%%%%%%%%%%%%%%%%%
\theorem \label{_vanishing_coho_weight_Theorem_}
Let $M$ be a compact Vaisman manifold, that is,
an LCK manifold with parallel Lee form $\theta$, and $V$ a positive
tensor power of the weight bundle (\ref{_weight_bu_Definition_}),
equipped with a Hermitian structure as in
Subsection \ref{_weight_curv_Subsection_}. Consider 
the $\omega_0$-differential 
$\bar\6_0:\; \Lambda^{0,p}(M)\otimes V \arrow \Lambda^{0,p+1}(M)\otimes V$
associated with the Chern connection in $V$ (see Subsection 
\ref{_Laplace_in_bundle_Subsection_}), and let
$\eta\in \Lambda^{0,p}(M)\otimes V$ be a 
non-zero $\bar\6_0$-harmonic form. Then $p\geq \dim_\C M -1$. 
Moreover, if $p= \dim_\C M -1$, then $\eta\cntrct \theta^\sharp =0$,
where $\theta^\sharp$ is the vector field dual to $\theta$.

\hfill

%%%%%%%%%%%%%%%%%%%%%%%%%%%%%%%%%%%%%%%%%%%%%%%%
\remark 
An elementary linear-algebraic argument implies that all
$(n-1,0)$-forms satisfying $\eta\cntrct \theta^\sharp =0$
are proportional to ${\goth L}\cntrct \theta^\sharp$,
where  ${\goth L}$ is a non-degenerate $(n,0)$-form on $M$.

\hfill

{\bf Proof of \ref{_vanishing_coho_weight_Theorem_}:}
By \ref{_Lapla_compa_Theorem_}, we have
\[
\Delta_{\6_0}-\Delta_{\bar\6_0} = \1 [ \Theta_{V}, \Lambda_0],
\]
where $\Theta_{V}$ is the curvature operator of $V$.
By \ref{_curva_L_omega_0_Theorem_},  
$\1\Theta_{V}= c \omega_0$, $c>0$.
Therefore, 
\begin{equation}\label{_Laplace_via_H_0_Equation_}
\Delta_{\bar \6_0} = \Delta_{\6_0} + c H_0,
\end{equation}
where $H_0 = [ L_0, \Lambda_0]$.
By \eqref{_H_0_explicitly_Equation_},
$H_0$ is positive definite 
on $(r,0)$-forms for $r< n-1$.
By \eqref{_Laplace_via_H_0_Equation_},
the Laplace operator $\Delta_{\bar\6_0}$
is a sum of a positive semidefinite
operator $\Delta_{\6_0}$ and
positive definite $c H_0$. Therefore,
all eigenvalues of $\Delta_{\bar \6_0}$
are strictly positive, and there are
no harmonic $(r,0)$-forms $r< n-1$. 

Similarly, $H_0$ is positive semidefinite
on $(n-1,0)$-forms, and its only zero eigenvalue
corresponds to the form 
\[ \nu= \xi_1\wedge \xi_2 \wedge ... \wedge \xi_{n-1}.
\]
(we use the notation of \ref{_omega_0_eigenvalues_Proposition_}
here). Therefore, if $\Delta_{\bar \6_0}(\eta)=0$,
then $\Delta_{\6_0}(\eta)=0$ and
$H_0(\eta)=0$, and 
$\eta$ is proportional to $\nu$.
We proved \ref{_vanishing_coho_weight_Theorem_}. \endproof

%%%%%%%%%%%%%%%%%%%%%%%%%%%%%%%%%%%%%%%%%%%%%%%%%%%%%%%%%%%%
\subsection{Basic cohomology}
%%%%%%%%%%%%%%%%%%%%%%%%%%%%%%%%%%%%%%%%%%%%%%%%%%%%%%%%%%%%

To be able to use the vanishing theorem
obtained above, we need a way
to compare the cohomology of $\bar\6_0$ and the Dolbeault 
cohomology. This comparison is obtained from the results of 
K. Tsukada (\cite{_Tsukada:forms_}).

\hfill

Let $M$ be a manifold equipped with a foliation
$\Xi$. A form $\eta$ is called {\bf basic} if
for all vector fields $v\in \Xi$, we have
$\eta\cntrct v=0, \Lie_v\eta=0$. Locally,
such forms are lifted from the space of leaves of
$\Xi$, if this space is defined. Clearly, $d\eta$
is basic if $\eta$ is basic. The de Rham
cohomology of basic forms is called 
{\bf the basic cohomology of the foliation
$\Xi$} (\cite{_Tondeur_}, \cite{_Nishikawa_Tondeur_}). 

\hfill

Given a holomorphic foliation on a complex manifold,
we cal also define the basic Dolbeault cohomology,
as the cohomology of the differential $\bar\6$ on
the basic forms. Clearly, on basic forms,
the differential $\bar\6$ is equal to the
differential $\bar\6_0$ defined above.

\hfill

%%%%%%%%%%%%%%%%%%%%%%%%%%%%%%%%%%%%%%%%%%%%%%%%%%%%%%%%%%%%
\theorem\label{_Tsuka_compa_harmo_and_basic_Theorem_}
Let $M$ be an $n$-dimensional compact Vaisman manifold,
and  $\eta$ a $(p,q)$-form, with $p+q\leq n-1$. 
Denote by $\theta^{0,1}$ the $(0,1)$-part of the
Lee form, and let $\Lambda_0$ be the Hodge operator
associated with the form $\omega_0$ (Subsection
\ref{_Koda_for_6_0_Subsection_}). Then
the following conditions are equivalent.

\begin{description}
\item[(i)] The form $\eta$ is $\bar\6$-harmonic: 
$\bar\6\bar\6^* \eta +  \bar\6^*\bar\6\eta=0$.
\item[(ii)] $\eta$ has a decomposition
$\eta=\theta^{0,1}\wedge \alpha +\beta$,
where $\alpha$, $\beta$ are basic forms
of the canonical foliation $\Xi$ which are
$\bar\6_0$-harmonic and satisfy 
$\Lambda_0 \alpha = \Lambda_0\beta=0$.
\end{description}

{\bf Proof:} This is \cite{_Tsukada:forms_}, Theorem 3.2. \endproof

%%%%%%%%%%%%%%%%%%%%%%%%%%%%%%%%%%%%%%%%%%%%%%%%%%%%%%%%%%%%
\subsection{The cohomology of the structure sheaf}
%%%%%%%%%%%%%%%%%%%%%%%%%%%%%%%%%%%%%%%%%%%%%%%%%%%%%%%%%%%%

The vanishing result of Subsection 
\ref{_vani_for_6_0_Subsection_}
can be used to prove the following theorem.

\hfill

%%%%%%%%%%%%%%%%%%%%%%%%%%%%%%%%%%%%%%%%%%%%%%%%
\theorem \label{_cohomo_calo_M_Theorem_}
Let $M$ be a compact Vaisman manifold, such that the
canonical bundle $K(M)$ is a negative power of the 
weight bundle\footnote{By \ref{_can_cla_weight_on_Einstein_Proposition_}, 
this holds for all Einstein-Weyl manifolds. Hence, any LCHK manifold
satisfies this assumption.}. Then the holomorphic 
cohomology of the structure sheaf satisfy
$H^i(\calo_M) =0$ for $i>1$, $\dim H^1(\calo_M) =1$.

\hfill

{\bf Proof:} \ref {_cohomo_calo_M_Theorem_}
is a consequence of 
\ref{_vanishing_coho_weight_Theorem_}
and \ref{_Tsuka_compa_harmo_and_basic_Theorem_}.
Indeed, \ref{_Tsuka_compa_harmo_and_basic_Theorem_} 
implies that to prove that $H^i(\calo_M) =0$ for $i>1$, 
$\dim H^1(\calo_M) =1$, it suffices to show that
all basic $\bar\6$-harmonic $(0,p)$-forms vanish, for
$p>0$. A basic form is $\bar\6$-harmonic
if and only if it is $\bar\6_0$-harmonic.
By Serre's duality, 
\[
{\cal H}_{\bar \6_0}^p(\calo_M) \cong 
\overline{{\cal H}_{\bar \6_0}^{n-p}(K)},
\]
(\ref{_Serre_s_for_6_0_Theorem_}).
By 
\ref{_vanishing_coho_weight_Theorem_},
${\cal H}_{\bar \6_0}^{n-p}(K)=0$ for $p>1$, and
for $p=1$ all non-trivial classes 
$\eta \in {\cal H}_{\bar \6_0}^{n-p}(K)$ satisfy
$\eta\cntrct {\theta^{0,1}}^\sharp =0$.
Applying the Serre's duality again,
we obtain that any $\bar\6_0$-harmonic
$(0,p)$-form $\eta_1\in {\cal H}_{\bar \6_0}^p(\calo_M)$ 
vanishes for $p>1$, and for $p=1$, 
$\eta_1$ is proportional to $\theta^{0,1}$.
Therefore, for $p=1$, $\eta$ cannot be basic,
and we obtain that there are no non-trivial
basic $\bar\6$-harmonic $(0,p)$-forms.
This proves \ref{_cohomo_calo_M_Theorem_}.
\endproof

\hfill

A similar theorem was obtained in \cite{_Aleksandrov_Ivanov_}
using an estimate of Ricci curvature.

%%%%%%%%%%%%%%%%%%%%%%%%%%%%%%%%%%%%%%%%%%%%%%%%%%%%%%%%%%%%

\section{Holomorphic forms on Einstein-Weyl LCK manifolds}

%%%%%%%%%%%%%%%%%%%%%%%%%%%%%%%%%%%%%%%%%%%%%%%%%%%%%%%%%%%%

In this Section, we prove that all
holomorphic $(p,0)$-forms on a compact Einstein-Weyl LCK
manifold vanish, for $p>0$.

%%%%%%%%%%%%%%%%%%%%%%%%%%%%%%%%%%%%%%%%%%%%%%%%%%%%%%%%%%%%
\subsection{The $\omega_0$-Yang-Mills bundles}
%%%%%%%%%%%%%%%%%%%%%%%%%%%%%%%%%%%%%%%%%%%%%%%%%%%%%%%%%%%%

Let $M$ be a Vaisman manifold, that is, 
an LCK manifold with parallel Lee form $\theta$, and
$\omega_0 = d^c(\theta)$ the standard 2-form
(Section \ref{_omega_0_Section_}).
For our purposes, $\omega_0$ plays
the role of a K\"ahler form on $M$.  It is 
natural to study the Yang-Mills geometry associated
with $\omega_0$.

\hfill

%%%%%%%%%%%%%%%%%%%%%%%%%%%%%%%%%%%%%%%%%%%%%%%%
\definition
In the above assumptions, let $B$ be a holomorphic Hermitian
bundle equipped with the standard (Chern) connection,
$\Theta_B$ its curvature, and 
$\Lambda_0:\; \Lambda^p(M) \arrow \Lambda^{p-2}(M)$
the Hodge-type operator associated with $\omega_0$
(see Subsection \ref{_Koda_for_6_0_Subsection_}).
We say that $B$ is {\bf $\omega_0$-Yang-Mills} if 
\begin{equation}\label{_omega_0_YM_Equation_}
\Lambda_0(\Theta_B) = -c\1\Id_B,
\end{equation}
where $Id_B$ is the identity section of $\End(B)$,
and $c$ a constant. We call $c$ {\bf the Yang-Mills
constant of $B$}. 

\hfill

The $\omega_0$-Yang-Mills bundles satisfy the same
elementary properties as the usual Yang-Mills bundles.
In particular, any tensor power of $\omega_0$-Yang-Mills bundles
is again $\omega_0$-Yang-Mills. The following theorem provides
us with an example of a $\omega_0$-Yang-Mills bundle.

\hfill

%%%%%%%%%%%%%%%%%%%%%%%%%%%%%%%%%%%%%%%%%%%%%%%%%%%%%%%%%%%%
\proposition\label{_TM_YM_for_Einstein_Weyl_Proposition_}
Let $M$ be an Einstein-Weyl LCK manifold 
with parallel Lee form.\footnote{
By \ref{_E-W_theta_par_Theorem_}, 
a compact Einstein-Weyl LCK manifold
with the Gauduchon metric 
has parallel Lee form.}
Consider the tangent bundle 
$TM$ as a holomorphic Hermitian bundle.
Then $TM$ is $\omega_0$-Yang-Mills,
with positive Yang-Mills constant.

\hfill

{\bf Proof:} Let $\nabla_W$ be the Weyl connection
on $TM$, and $\nabla_C$ the Chern connection.
By definition, $\nabla^{0,1}_W = \nabla^{0,1}_C$.
Consider the weight bundle $(L, \nabla_L)$ on $M$. 
Since $TM$ is equipped with a metric with values in
$L^{\otimes 2}$, the natural connection
$\nabla_{W,L}$ on $TM\otimes (L^*)^{\otimes 2}$ induced by $\nabla_W$, 
$\nabla_L$, is Hermitian. Therefore,
it is a Chern connection on 
the holomorphic Hermitian
vector bundle $TM\otimes L^{\otimes 2}$.

Since $\nabla_L$ is flat, we have
$\nabla_{W,L}^2 = \nabla_{W}^2$. The Chern connection
on $TM\otimes (L^*)^{\otimes 2}$ is obtained as a tensor product
of the Chern connection on $TM$ and that on $L$.
We obtain
\[ \nabla_{W}^2=\nabla_{W,L}^2 = \nabla_{C}^2 - 2 C,
\]
where $C = -2\1 \omega_0$ is the curvature of the Chern
connection on $L$ (see \ref{_curva_L_omega_0_Theorem_}).
This gives
\[ 
\nabla_{C}^2 = 2 C + \nabla_{W}^2.
\]
Since $\Lambda_0(C)$ is a constant,
to prove \ref{_TM_YM_for_Einstein_Weyl_Proposition_}
it remains to show that
\begin{equation}\label{_lambda_from_curv_nabla_W_Equation_}
\Lambda_0(\Theta_{W}) =0,
\end{equation}
where $\Theta_W \in \Lambda^2(M) \otimes \End(TM)$
is the curvature of the Weyl connection.

\hfill

The following lemma is clear.

\hfill

%%%%%%%%%%%%%%%%%%%%%%%%%%%%%%%%%%%%%%%%%%%%%%%%
\lemma\label{_curva_omega_primiti_Lemma_}
Let $M$ be an Einstein-Weyl LCK manifold
equipped with a Gauduchon metric
and $\omega$ the Hermitian form on $M$.
Denote by $\Lambda_\omega:\; \Lambda^p(M) \arrow \Lambda^{p-2}(M)$
the Hermitian adjoint  operator to $\eta\arrow \eta\wedge \omega$.
Then $\Lambda_\omega(\Theta_{W}) =0$.

\hfill

{\bf Proof:} Let $\tilde M$ be the 
universal covering of $M$, and $rg$ the Calabi-Yau metric on
$\tilde M$. 
The form $\Theta_{W}$ is the curvature
of the Levi-Civita connection of the K\"ahler
metric on $\tilde M$.  Therefore, $TM$ is Yang-Mills
with respect to this metric, and $\Theta_{W}$ is orthogonal to 
the K\"ahler form $r\omega$. This implies that 
$\Theta_{W}$ is orthogonal to $\omega$. 
We obtain $\Lambda_\omega(\Theta_{W}) =0$.
\endproof

\hfill

Return to the proof of 
\ref{_TM_YM_for_Einstein_Weyl_Proposition_}.
By \ref{_curva_omega_primiti_Lemma_}, we have
\[ \Lambda_\omega(\Theta_{W}) =0,\] and we need to show
\[\Lambda_0(\Theta_{W}) =0.\] 
By \ref{_omega_0_eigenvalues_Proposition_},
we have
\[ \Lambda_{\omega}\Theta_W - \Lambda_0\Theta_W =
   \Theta_W(\theta^\sharp, I(\theta^\sharp)).
\]
Therefore, to prove 
$\Lambda_{\omega}\Theta_W =\Lambda_0\Theta_W=0$, we need
to show that the curvature $\Theta_W$ restricted to 
the canonical foliation $\Xi$ vanishes. 
The following lemma finishes the proof of
\ref{_TM_YM_for_Einstein_Weyl_Proposition_}.

\hfill

%%%%%%%%%%%%%%%%%%%%%%%%%%%%%%%%%%%%%%%%%%%%%%%%
\lemma \label{_conne_flat_on_Xi_Lemma_}
Let $M$ be a Vaisman manifold,
$\Xi\subset TM$ the canonical foliation,
and $\nabla$ the Weyl connection. Then
$\nabla$ is flat on the leaves of 
$\Xi$.

\hfill

{\bf Proof:} Using 
\eqref{_nabla_Weyl_via_nabla_g_Equation_} and $\nabla_g\theta^\sharp=0$, 
we obtain
\begin{equation}\label{_Weyl_conne_to_Theta_Equation_}
\nabla_X \theta^\sharp = X,
\end{equation}
for all vector fields $X\in TM$.
Therefore,
\[
\nabla_X \nabla_Y\theta^\sharp = \nabla_X Y - \nabla_Y X - [X, Y] =0.
\]
We find that $\Theta_W(X, Y, \theta^\sharp)=0$
for all $X, Y\in TM$. Using the symmetries of a curvature tensor,
we obtain that $\Theta_W(X, \theta^\sharp, Y)=0$
identically. This proves  
\ref{_conne_flat_on_Xi_Lemma_}. 
\ref{_TM_YM_for_Einstein_Weyl_Proposition_} is proven. \endproof

%%%%%%%%%%%%%%%%%%%%%%%%%%%%%%%%%%%%%%%%%%%%%%%%%%%%%%%%%%%%%%%%%%%%%%%%
\subsection{Vanishing theorem for negative bundles}
%%%%%%%%%%%%%%%%%%%%%%%%%%%%%%%%%%%%%%%%%%%%%%%%%%%%%%%%%%%%%%%%%%%%%%%%

In this Subsection we present another vanishing theorem based
on \ref{_Lapla_compa_Theorem_}.
More general results are possible, in line with the
standard vanishing theorems in algebraic geometry.

\hfill

%%%%%%%%%%%%%%%%%%%%%%%%%%%%%%%%%%%%%%%%%%%%%%%%%%%%%%%%%%%%
\theorem\label{_nega_bu_no_sec_Theorem_}
Let $M$ be a compact Vaisman manifold, $B$ a holomorphic Hermitian
bundle and $\Theta_B \in \Lambda^{1,1}(M) \otimes \End(B)$
its curvature. Consider the $\omega_0$-Hodge operator
\[ \Lambda_0:\; \Lambda^{1,1}(M) \otimes \End(B)\arrow\End(B)
\]
(Subsection \ref{_Koda_for_6_0_Subsection_}).
Assume that the self-adjoint operator $\1 \Lambda_0(\Theta_B)$
is strictly negative everywhere in $M$.
Then $B$ has no non-zero holomorphic sections.

\hfill

{\bf Proof:} 
Let $\beta\in B$ be a holomorphic section of $B$.
Then $\nabla^{0,1}(\beta)=0$. From the definition
of $\bar\6_0$, we obtain that $\bar\6_0\beta=0$.
Clearly, $\bar\6_0^*\beta=0$ as well.
Therefore, $\beta$ is $\bar\6_0$-harmonic.
{}From the equality \eqref{_Lapla_compa_Equation_},
we find immediately that
\[
\bigg( \1 \Lambda_0 \Theta_B (\beta), \beta\bigg)> 0.
\]
This is impossible, because the operator 
$\1 \Lambda_0(\Theta_B)$ is strictly negative. \endproof

\hfill

%%%%%%%%%%%%%%%%%%%%%%%%%%%%%%%%%%%%%%%%%%%%%%%%%%%%%%%%%%%%
\corollary\label{_YM_bu_no_sections_Corollary_}
Let $M$ be a compact Vaisman manifold, not Kaehler,
and $B$ a $\omega_0$-Yang-Mills bundle with negative Yang-Mills
constant. Then $B$ has no global holomorphic sections.

\hfill

{\bf Proof:} Follows immediately from
\ref{_nega_bu_no_sec_Theorem_}.
\endproof

%%%%%%%%%%%%%%%%%%%%%%%%%%%%%%%%%%%%%%%%%%%%%%%%%%%%%%%%%%%%
\subsection{Vanishing for holomorphic forms} 
%%%%%%%%%%%%%%%%%%%%%%%%%%%%%%%%%%%%%%%%%%%%%%%%%%%%%%%%%%%%

The main result of this section is the following theorem.

\hfill

%%%%%%%%%%%%%%%%%%%%%%%%%%%%%%%%%%%%%%%%%%%%%%%%
\theorem \label{_holo_forms_vanish_Theorem_}
Let $M$ be a compact Einstein-Weyl LCK manifold.
Assume that $M$ is not K\"ahler.
Then all holomorphic $p$-forms 
on $M$ vanish, for all $p>0$. 

\hfill

{\bf Proof:} Follows 
{}from \ref{_TM_YM_for_Einstein_Weyl_Proposition_}
and \ref{_YM_bu_no_sections_Corollary_}.\endproof

\hfill

A similar theorem was obtained in \cite{_Aleksandrov_Ivanov_}
using an estimate of Ricci curvature.

\hfill

This result has the following topological implications.

\hfill

%%%%%%%%%%%%%%%%%%%%%%%%%%%%%%%%%%%%%%%%%%%%%%%%%%%%%%%%%%%%
\theorem\label{_h^1=1_Theorem_}
Let $M$  be a compact Einstein-Weyl LCK manifold.
Assume that $M$ is not K\"ahler. Then the first Betti
number of $M$ is 1: \[ h^1(M)=1.\]

{\bf Proof:} Consider the Dolbeault spectral sequence
$E^{p,q}_r$ associated with $M$. Then
$E^{p,q}_2 = H^q(\Omega^p(M))$.
By \ref{_holo_forms_vanish_Theorem_}, we have
$H^0(\Omega^1(M))=0$, by \ref{_cohomo_calo_M_Theorem_}, we have 
$H^1(\Omega^0(M))=\C$. Since the Dolbeault spectral sequence
converges to de Rham cohomology, we obtain that
\[ h^1(M) \leq \dim H^0(\Omega^1(M)) + \dim H^1(\Omega^0(M)) =1.\]
It remains to prove that $h^1(M)\neq 0$.
The monodromy of the weight bundle lies in $\R^{>0}$,
hence it is abelian and torsion-free.
If $h^1(M)=0$, the weight bundle has 
trivial monodromy. By \ref{_theta_exact_confo_Riemann_Claim_},
an LCK manifold with trivial weight bundle
is K\"ahler. Since $M$ is not Kaehler, we have
$h^1(M)>0$. This proves \ref{_h^1=1_Theorem_}.
\endproof

\hfill

A similar theorem was obtained in \cite{_Aleksandrov_Ivanov_}
using an estimate of Ricci curvature.

%%%%%%%%%%%%%%%%%%%%%%%%%%%%%%%%%%%%%%%%%%%%%%%%%%%%%%%%%%%%

\section{Sasakian geometry: an introduction}

%%%%%%%%%%%%%%%%%%%%%%%%%%%%%%%%%%%%%%%%%%%%%%%%%%%%%%%%%%%%

Sasakian manifolds were introduced by S. Sasaki
(\cite{_Sasaki_}). In this Section we reproduce the
definition and some basic results on 1-Sasakian and 3-Sasakian
manifolds. For a survey and reference of Sasakian geometry,
the reader should consult \cite{_Boyer_Galicki_}

\hfill

%%%%%%%%%%%%%%%%%%%%%%%%%%%%%%%%%%%%%%%%%%%%%%%%
\definition
Let $(X,g)$ be a Riemannian manifold.
A {\bf cone} of $X$ is a Riemannian manifold
$C(X):= (X \times \R^{>0}, dt^2 + t^2 g )$,
where $t$ is the parameter on $\R^{>0}$. 
For any $\lambda\in \R^{>0}$,
the map 
\[ \tau_\lambda:\; C(X) \arrow C(X), \ \  (x, t) \arrow (x, \lambda t)
\]
multiplies the metric by $\lambda^2$.

\hfill

%%%%%%%%%%%%%%%%%%%%%%%%%%%%%%%%%%%%%%%%%%%%%%%%
\definition
Let $X$ be a Riemannian manifold. A 1-Sasakian structure in $X$
is a complex structure on $C(X)$ satisfying the following
\begin{description}
\item[(i)] The metric on $C(X)$ is K\"ahler
\item[(ii)] The map $\tau_\lambda:\; C(X) \arrow C(X)$ is holomorphic,
for all $\lambda\in \R^{>0}$.
\end{description}
A 3-Sasakian structure on $X$ is a hypercomplex strucure on
$C(X)$ such that
\begin{description}
\item[(i)] The metric on $C(X)$ is hyperk\"ahler
\item[(ii)] The map $\tau_\lambda:\; C(X) \arrow C(X)$ is 
compatible with the hypercomplex structure,
for all $\lambda\in \R^{>0}$.
\end{description}

\hfill

%%%%%%%%%%%%%%%%%%%%%%%%%%%%%%%%%%%%%%%%%%%%%%%%
\remark
A 3-Sasakian manifold is equipped with a 1-Sasakian structure,
for any quaternion $L\in {\Bbb H}, L^2 =-1$.

\hfill

One can define 1-Sasakian and 3-Sasakian structures
in terms of vector fields $I(dt^\sharp)\in TX$ and
$I(dt^\sharp)$, $J(dt^\sharp)$, $K(dt^\sharp)\in TX$,
and the associated contact structures (see \cite{_Boyer_Galicki_}).
Historically, these manifolds were 
defined in this fashion. The Sasakian geometry 
relates to contact geometry in exactly
the same way as K\"ahler geometry is related to 
symplectic geometry. 

\hfill

%%%%%%%%%%%%%%%%%%%%%%%%%%%%%%%%%%%%%%%%%%%%%%%%
\definition 
Let $X$ be a 1-Sasakian manifold, such that the Riemannian
metric on $X$ satisfies Einstein equation. Then 
$X$ is called {\bf Sasa\-kian-\-Ein\-stein}. 

\hfill

Sasakian-Einstein manifolds can be characterized in terms
of their cones, as follows. 

\hfill

%%%%%%%%%%%%%%%%%%%%%%%%%%%%%%%%%%%%%%%%%%%%%%%%
\proposition \label{_Sasa_E-Calabi-Yau_cone_Proposition_}
Let $X$ be a 1-Sasakian manifold, $\dim_\R X = 2n+1$.
Then $X$ is Sasakian-Einstein if and only if its cone is Ricci-flat.
In this case the Einstein constant of $X$ is equal to $2n$.

\hfill

{\bf Proof:} See e.g. \cite{_Boyer_Galicki_}, Proposition 1.1.9.
\endproof

\hfill

%%%%%%%%%%%%%%%%%%%%%%%%%%%%%%%%%%%%%%%%%%%%%%%%
\remark A cone of a 3-Sasakian manifold is hyperk\"ahler,
hence Ricci-flat. Therefore, a 3-Sasakian 
manifold is always Einstein.

\hfill

%%%%%%%%%%%%%%%%%%%%%%%%%%%%%%%%%%%%%%%%%%%%%%%%
\remark\label{_Myers_Remark_}
By Myers' theorem, a complete Einstein manifold with
positive Einstein constant is compact and has finite
fundamental group. 

%%%%%%%%%%%%%%%%%%%%%%%%%%%%%%%%%%%%%%%%%%%%%%%%

\section{Vaisman manifolds with $h^1(M)=1$ 
and Sa\-sa\-kian geometry}

%%%%%%%%%%%%%%%%%%%%%%%%%%%%%%%%%%%%%%%%%%%%%%%%

%%%%%%%%%%%%%%%%%%%%%%%%%%%%%%%%%%%%%%%%%%%%%%%%
\subsection{Vaisman manifolds with exact Lee form}
\label{_exa_Lee_Sasakian_Subsection_}
%%%%%%%%%%%%%%%%%%%%%%%%%%%%%%%%%%%%%%%%%%%%%%%%

Vaisman manifolds are intimately related to Sasakian geometry,
as the following proposition shows.

\hfill

%%%%%%%%%%%%%%%%%%%%%%%%%%%%%%%%%%%%%%%%%%%%%%%%%%%%%%%%%%%%%%%%%%%%%%%%
\proposition\label{_Vaisman_decompo_Sasa_Proposition_}
(\cite{_Kamishima_Ornea_}, \cite{_Gini_Ornea_Parton_})
Let $M$ be a Vaisman manifold, that is, an LCK manifold
with a parallel Lee form $\theta$. Assume that $\theta\neq 0$
and $\theta$ is exact: $\theta = d\zeta$. Assume, moreover, 
that the Gauduchon metric on $M$ is complete.
Consider the function $\zeta$
as a map $\zeta:\; M \arrow \R$.
Then $M$ is isometric to a product
$X \times \R$, where $X$ a complete 
1-Sasakian manifold, with the projection
to $M = X \times \R\arrow \R$ given
by $\zeta$. Moreover, if 
$M$ is LCHK, then $X$ is naturally
a 3-Sasakian manifold.

\hfill

{\bf Proof:} 
Consider the vector field $\theta^\sharp$
dual to $\theta$. By definition,
$\theta^\sharp$ is a parallel vector field on $M$,
and  $\theta^\sharp$ is equal to the gradient of $\zeta$.
Therefore, the gradient flow $\Phi_\lambda$ associated with $\zeta$ is 
an isometry. The map $\Phi_\lambda$  is well defined
because $M$ is complete.
Beind a gradient flow, $\Phi_\lambda$ commutes with $\zeta$:
\begin{equation}\label{_flow_and_phi_Equation_}
\zeta(\Phi_\lambda(x)) = \lambda+ \zeta(x). 
\end{equation}
Therefore,
$\Phi_\lambda$ induces an isometry on the fibers of 
$\zeta(x)$. Denote by $X$ the fiber $\zeta^{-1}(0)$.
Let $R:\; X \times \R \arrow M$ map
$(t, x)$ into $\Phi_t(x)$. 
Clearly, $R$ is an isometry. Denote the Riemannian metric on $M$ 
by $g$. By \ref{_kah_form_on_LCK_w_pote_Claim_},
the metric $e^{2t} g$ is 
K\"ahler. After a reparametrization
$t \arrow e^t$, we obtain that the metric
$(e^{2t} g, (d e^t)^2)$ is the cone metric on
$M = X \times \R^{>0}$. Therefore,
$X$ is 1-Sasakian. If $M$ is LCHK,
then $C(X)$ is hyperk\"ahler 
(\ref{_LCHK_w_trivial_monodro_Claim_}),
hence $X$ is 3-Sasakian.
\endproof

\hfill

%%%%%%%%%%%%%%%%%%%%%%%%%%%%%%%%%%%%
\remark
A covering of a compact manifold is complete.
Therefore, \ref{_Vaisman_decompo_Sasa_Proposition_}
holds for a covering $\tilde M$ of any 
compact Vaisman manifold, if $\tilde M$ has
exact Lee form.

%%%%%%%%%%%%%%%%%%%%%%%%%%%%%%%%%%%%%%%%%%%%%%%%
\subsection{Structure theorem for Vaisman manifolds with 
the first Betti number 1}
%%%%%%%%%%%%%%%%%%%%%%%%%%%%%%%%%%%%%%%%%%%%%%%%

It is possible to classify the Vaisman manifolds with 
the first Betti number 1, in terms of 1-Sasakian geometry, as follows.

\hfill

Let $X$ be a 1-Sasakian manifold, and $C(X)$ its cone. Given a number
 $q\in \R$, $q>1$, consider an equivalence relation
$\sim_q$ on $C(X) = X\times \R^{>0}$ generated
by $(x,t) \sim (x, qt)$. Since the
map $(x,t) \arrow (x, qt)$ mupltiplies the metric
by $q^2$, the quotient $M = C(X)/\sim_q$
is an LCK manifold. Moreover,
$M$ is a Vaisman manifold, with the Gauduchon
metric provided by an isomorphism
$M \cong X \times S^1$.

This construction can be generalized as follows.
Let $\phi:\; X \arrow X$ be an automorphism of
a Sasakian structure. The map $(x, t) \stackrel{\phi_q} \arrow (\phi(x), qt)$
is compatible with the complex structure and multiplies
the metric by $q^2$. Therefore, the quotient
$M_{\phi, q}$ of $C(X)$ by the
corresponding equivalence relation
$\sim_{\phi,q}$ is an LCK manifold.
The following theorem is quite elementary.

\hfill

%%%%%%%%%%%%%%%%%%%%%%%%%%%%%%%%%%%%%%%%%%%%%%%%
\theorem\label{_Vaisman_stru_H^1=1_Theorem_}
(Structure theorem for Vaisman manifolds
with the first Betti number 1).
Let $X$ be a compact 1-Sasakian manifold,
$\phi:\; X \arrow X$ a 1-Sasakian automorphism,
and $M_{\phi, q}$  an LCK manifold constructed above.
Then $M_{\phi, q}$ is a Vaisman manifold satisfying 
the following conditions.
\begin{description}
\item[(i)] $h^1(M_{\phi, q})=1 \Leftrightarrow h^1(X)=0$
\item[(ii)]  
The 1-Sasakian manifold $X$,
together with the automorphism $\phi$, is uniquely
(up to a scaling) determined by the LCK structure on 
$M_{\phi, q}$.

\item[(iii)] Any compact Vaisman manifold
$M$, $h^1(M)=1$ can be constructed this way.
\end{description}

%%%%%%%%%%%%%%%%%%%%%%%%%%%%%%%%%%%%
\remark 
In \cite{_Gini_Ornea_Parton_}, it was shown that
$M_{\phi, q}$ is a Vaisman manifold (Proposition 7.4).

\hfill

{\bf Proof of \ref{_Vaisman_stru_H^1=1_Theorem_}:} 
To show that $M_{\phi, q}$ 
is a Vaisman manifold, consider the 
map $Id \times \log:\; C(X) \arrow X \times \R$,
$(x,t) \arrow (x, \log t)$. Let 
$\tilde g$ be the product metric on $X \times \R$ 
pulled back to $C(X)$.
The map $\phi_q$ induces an isometry,
hence $\tilde g$ corresponds to 
a metric $g$ on $M_{\phi, q}$.
By construction, $g$ belongs
to the same conformal class as
the LCK structure on $M_{\phi, q}$.
An elementary computation 
shows that $\nabla(\tilde g) = g\otimes dt$,
where $t$ is the parameter on 
$X \times \R$ corresponding to 
the second component. Therefore,
the Lee form of $g$ is parallel.

To find $H^1(M_{\phi, q})$, consider the 
natural projection $r:\; C(X)\arrow \R^{>0}$,
and let $\zeta = \log r$. The map
$\zeta$ sends  the points \[ (x, t)\sim_{\phi, q}(x', t')\]
to \[ \log t, \log t + k\log q, \ \ k\in \Z.\] Therefore,
$\zeta$ induces a map
$\zeta_q:\; M_{\phi, q}\arrow \R/(\log q)\Z$
{}from $M_{\phi, q}$ to a circle. By
construction, $d\zeta_q=\theta$, where
$\theta$ is the Lee form. Therefore,
\begin{equation}\label{_phi_q_fibra_Equation_}
\zeta_q:\; M_{\phi, q}\arrow S^1
\end{equation} 
is a smooth fibration, with the fiber $X$. 

Consider the Serre's spectral sequence 
$E^{p,q}_r$ for the fibration \eqref{_phi_q_fibra_Equation_}.
Since $H^i(S^1)=0$ for $i>0$, the spectral sequence 
$E^{p,q}_r$ degenerates in $E_2$. Therefore,
\begin{equation}\label{_fibra_over_S^1_cohomo_Equation_}
h^1(M_{\phi,q})) = h^1(S^1) + h^1(X) = h^1(X)+1
\end{equation}
This proves \ref{_Vaisman_stru_H^1=1_Theorem_} (i). 

To recover $(X, \phi)$ from the LCK structure on
$M= M_{\phi, q}$, notice that the Gauduchon
metric on $M$ is unique, hence the form $\theta$
is determined uniquely from the LCK geometry.
Applying \ref{_Vaisman_decompo_Sasa_Proposition_}
to $\tilde M = C(X)$, we reconstruct the
1-Sasakian structure on $X$, together
with an isomorpism $\tilde M \cong C(X)$.
The deck transform of $\tilde M$ induces
an automorphism $\phi$ of $X$. This
allows one to recover $(X, \phi)$
{}from the LCK structure on $M$.

It remains to prove \ref{_Vaisman_stru_H^1=1_Theorem_} (iii).

Let $M$ be a compact Vaisman manifold, $h^1(M)=1$,
$L$ its weight bundle and $\tilde M$ the covering
associated with the monodromy group $G$ of $L$.
Since $h^1(M)=1$, $M$ is non-K\"ahler, hence 
$L$ is non-trivial.
The monodromy group $G$ is naturally a 
subgroup of $\R^{>0}$. Therefore
$G$ is torsion-free and abelian. 
Since  $h^1(M)=1$, $G=\Z$.

We find that $M = \tilde M /\Z$. 
Let $\zeta:\; \tilde M \arrow \R$ be the function 
satisfying $d\zeta =\theta$. 
By \ref{_Vaisman_decompo_Sasa_Proposition_},
$\tilde M = X \times \R$, where $X$ is a complete 
1-Sasakian manifold, and the projection to the
second component given by $\zeta$.
Fix a point $x_0\in \tilde M$.
 Given $y\in \tilde M$
and a path $\gamma$ from $x_0$ to $y$,
we have 
\begin{equation}\label{_phi_integral_Equation_}
\zeta(x_0)-\zeta(y) = \int_\gamma \theta
\end{equation}
by the Stoke's formula.
Let $\gamma_0$ be the generator of 
the monodromy group of $L$ and 
$w$ the integral $w:=\int_{\gamma_0} \theta$.
Denote by $R:\; \tilde M \arrow \tilde M$ be the
monodromy transfrorm of $\tilde M$ associated with
$\gamma_0 \in G$.
Using an isomorphism $\tilde M = X \times \R$
and \eqref{_phi_integral_Equation_}, we obtain that
$R$ maps $(x, t)$ to $(x_1, t+w)$. This gives a map
$\phi:\; X \arrow X$, $x \arrow x_1$.
Clearly, $M= C(X)/\sim_{\phi, e^w}$.
 This proves \ref{_Vaisman_stru_H^1=1_Theorem_} (iii). \endproof

%%%%%%%%%%%%%%%%%%%%%%%%%%%%%%%%%%%%%%%%%%%%%%%%%%%%%%%%%%%%
\subsection{Structure theorem for Einstein-Weyl LCK and 
LCHK manifolds}
%%%%%%%%%%%%%%%%%%%%%%%%%%%%%%%%%%%%%%%%%%%%%%%%%%%%%%%%%%%%

Comparing \ref{_h^1=1_Theorem_} 
and \ref{_Vaisman_stru_H^1=1_Theorem_},
we immediately obtain the following structure
theorems.

\hfill

%%%%%%%%%%%%%%%%%%%%%%%%%%%%%%%%%%%%%%%%%%%%%%%%
\theorem 
Let $M$ be a compact Einstein-Weyl LCK 
manifold which is not K\"ahler. Then there exists
an Einstein Sasakian manifold $X$ and a
Sasakian automorphism $\phi:\; X \arrow X$
such that $M \cong C(X)/\sim_{\phi, q}$, 
for some $q\in \R, q >1$,
where $\sim_{\phi, q}$ is an equivalence relation
generated by $(x,t)\sim_{\phi, q} (\phi(x), qt)$.
Moreover, the manifold $X$ and an automorphism
$\phi$ are determined uniquely from the
LCK structure on $M$, up to a rescaling
of a  metric on $X$.

\hfill

{\bf Proof:}  By 
\ref{_Vaisman_decompo_Sasa_Proposition_}
the covering $\tilde M$ is isomorphic to $C(X)$,
for a 1-Sasakian manifold $X$. By \ref{_Calabi_yau_covering_Einstein_},
$C(X)$ is equipped with a Ricci-flat
K\"ahler metric. By \ref{_Sasa_E-Calabi-Yau_cone_Proposition_},
$X$ is Sasakian-Einstein. 
By \ref{_h^1=1_Theorem_}, $h^1(M)=1$.
By \ref{_Vaisman_stru_H^1=1_Theorem_},
$M \cong C(X)/\sim_{\phi, q}$ and $X, \phi$ are 
determined uniquely. 
\endproof

\hfill

%%%%%%%%%%%%%%%%%%%%%%%%%%%%%%%%%%%%%%%%%%%%%%%%
\theorem \label{_LCHK_Structure_Theorem_}
Let $M$ be a compact  LCHK 
manifold which is not hyperk\"ahler. Then there exists
a 3-Sasakian manifold $X$ and a
3-Sasakian automorphism $\phi:\; X \arrow X$
such that $M \cong C(X)/\sim_{\phi, q}$,
for some $q\in \R, q >1$,
where $\sim_{\phi, q}$ is an equivalence relation
generated by $(x,t)\sim_{\phi, q} (\phi(x), qt)$.
Moreover, the manifold $X$ and an automorphism
$\phi$ are determined uniquely from the
LCHK structure on $M$, up to a rescaling
of a  metric on $X$.

\hfill

{\bf Proof:} Since $M$ is LCHK, 
it is an Einstein-Weyl manifold. 
By \ref{_h^1=1_Theorem_}, $h^1(M)=1$.
By \ref{_Vaisman_stru_H^1=1_Theorem_},
$M \cong C(X)/\sim_{\phi, q}$ and $X, \phi$ are 
determined uniquely. The covering 
$\tilde M\cong C(X) $ is by definition hyperk\"ahler,
hence $X$ is 3-Sasakian. \endproof

%%%%%%%%%%%%%%%%%%%%%%%%%%%%%%%%%%%%%%%%%%%%%%%%
\subsection{Quasiregular LCHK manifolds}
%%%%%%%%%%%%%%%%%%%%%%%%%%%%%%%%%%%%%%%%%%%%%%%%

%%%%%%%%%%%%%%%%%%%%%%%%%%%%%%%%%%%%%%%%%%%%%%%%
\definition
Let $M$ be a Vaisman manifold, and $\theta^\sharp$
the vector field dual to the Lee form $\theta$.
Consider the flow $\Phi_t$ associated with $\theta^\sharp$.
The manifold $M$ is called {\bf quasiregular}
if for all compact sets $K \subset M$ and all points
$x\in M$, the intersection of the orbit
\[ V_x = \{ \Phi_t(x), t\in \R\}
\]
with $K$ is compact. In other words,
$M$ is quasiregular if the set $V_x$ does not have
concentration points outside itself, for all $x\in M$.

\hfill

Given an LCHK manifold $M$, one can
construct a number of foliations on $M$, similar
to the canonical foliation $\Xi$.
The leaf space of these foliations
will be an orbifold if $M$ is 
quasiregular. This way, we can reduce
a quasiregular LCHK manifold to
\begin{description}
\item[(a)] A 3-Sasakian orbifold (by taking a 
leaf space of the real 1-dimensional foliation generated
by $\theta^\sharp$).
\item[(b)] An holomorphic contact K\"ahler-Einstein 
orbifold (by taking a leaf space of $\Xi$)
\item[(c)] A quaternionic K\"ahler orbifold
(by taking a 
leaf space of the real 4-dimensional foliation generated
by $\theta^\sharp$, $I(\theta^\sharp)$, 
$J(\theta^\sharp)$ $K(\theta^\sharp)$).
\end{description}
For details of these constructions
and further results see \cite{_Ornea:LCHK_}, 
\cite{_Ornea_Piccini_}. 

\hfill

Using the structure theorem 
(\ref{_LCHK_Structure_Theorem_}),
it is possible to determine the class of
quasiregular LCHK manifolds in terms of 3-Sasakian
fibrations.
The following claim is clear.

\hfill

%%%%%%%%%%%%%%%%%%%%%%%%%%%%%%%%%%%%
\claim
Let $M$ be a compact LCHK manifold,
obtained as in \ref{_LCHK_Structure_Theorem_}
{}from a 3-Sasakian manifold $X$ and a 3-Sasakian
automorphism $\phi:\; X \arrow X$. Then 
$M$ is quasiregular if and only if 
$\phi$ is a finite order automorphism.

\endproof

%%%%%%%%%%%%%%%%%%%%%%%%%%%%%%%%%%%%%%%%%%%%%%%%

\section{Appendix A. $h^1(M)=1$ for 
Ein\-stein-\-Weyl LCK manifolds}
\label{_h^1=1_Appendix_}

%%%%%%%%%%%%%%%%%%%%%%%%%%%%%%%%%%%%%%%%%%%%%%%%

In this Appendix, we give a direct proof of a 
version of \ref{_h^1=1_Theorem_}.

\hfill

Let $M$ be a complete Vaisman manifold, and $\tilde M$ the covering 
associated with the monodromy of the weight bundle. 
Then $\tilde M = C(X)$, for a complete 1-Sasakian manifold
$X$ (\ref{_Vaisman_decompo_Sasa_Proposition_}). 
Whenever $M$ is Einstein-Weyl, the manifold $X$ becomes
Sasakian-Einstein  (this is implied immediately by 
\ref{_Calabi_yau_covering_Einstein_}, 
\ref{_Sasa_E-Calabi-Yau_cone_Proposition_}). 
By Myers' Theorem (\ref{_Myers_Remark_}), then, $X$ is compact,
and $\pi_1(X)$ is finite. Therefore,
the following theorem implies \ref{_h^1=1_Theorem_}.

\hfill

%%%%%%%%%%%%%%%%%%%%%%%%%%%%%%%%%%%%%%%%%%%%%%%%
\theorem\label{_h^1=1_seco_proof_Theorem_}
Let $M$ be a compact Vaisman manifold, and 
$\tilde M$ its covering associated with the monodromy $G$ of the
weight bundle. We have $\tilde M \cong C(X)$,
for some 1-Sasakian manifold $X$ (\ref{_Vaisman_decompo_Sasa_Proposition_}).
Assume that $X$ is compact and $h^1(X)=0$. Then $h^1(M)=1$.

\hfill

{\bf Proof:} By definition, $M = \tilde M/G$. 
Let $\zeta:\; \tilde M \arrow \R$ 
be a function such that $d\zeta=\theta$.
Denote by $\chi:\; G \arrow \R$ is the group homomorphism
$\gamma \arrow \int_\gamma \theta$. 
It is easy to check that this map 
is a logarithm of the monodromy map $G \arrow \R^{>0}$.
By \eqref{_phi_integral_Equation_}, the monodromy
acts in such a way that for all $\gamma\in G$
we have 
\begin{equation}\label{_zeta_and_chi_homo_Equation_}
   \zeta (\gamma x) = \chi(\gamma) + \zeta(x).
\end{equation}
Since $G$ is the monodromy of $L$, 
the map $\chi$ is a monomorphism. 
Now, either $G=\Z$, and in this case
$M$ is fibered over a circle with fibers $X$, hence
$h^1(M) = h^1(X)+1 = 1$ 
(see \eqref{_fibra_over_S^1_cohomo_Equation_}); 
or $\chi(G)$ is dense in $\R$. 
To prove \ref{_h^1=1_seco_proof_Theorem_}
it remains to show that
$\chi(G)$ cannot be dense. 

Consider an interval $[0, 1] \subset \R$, and let
$\tilde M_0:= \zeta^{-1}([0,1])$
be the corresponding subset in $\tilde M$. Clearly
$\tilde M_0= X \times [0,1]$ is compact.
Given a point $x \in \tilde M$, $\zeta(x)=0$, let $Gx$ denote 
its orbit with respect to the monodromy action. 
By \eqref{_zeta_and_chi_homo_Equation_},
$Gx$ meets $\tilde M_0$ for all $\gamma\in G$
such that $\chi(\gamma) \in [0,1]$. If 
$\chi(G)$ is dense, this set is infinite.
We obtain that $Gx \cap \tilde M_0$ is infinite.
Since this set is compact, $Gx$ has concentration points.
This is clearly impossible, because $\tilde M \arrow \tilde M/G$
is a covering. Therefore, $\chi(G)$ cannot be dense.
We proved \ref{_h^1=1_seco_proof_Theorem_}.
\endproof

%%%%%%%%%%%%%%%%%%%%%%%%%%%%%%%%%%%%%%%%%%%%%%%%

\section{Appendix B. Counterexamples}

%%%%%%%%%%%%%%%%%%%%%%%%%%%%%%%%%%%%%%%%%%%%%%%%

For a general Vaisman manifold (without the
Einstein-Weyl assumption), \ref{_holo_forms_vanish_Theorem_} 
and \ref{_h^1=1_Theorem_} are 
false, as the following example shows. 

Let $S$ be an elliptic curve and $B$ a negative holomorphic line bundle.
Assume that $B$ is equipped with a Hermitian metric in such 
a way that the corresponding Chern connection has curvature
\begin{equation}\label{_Theta_B_Kah_Equation_}
\Theta_B = c \1 \omega,
\end{equation}
where $c >0$ is a positive constant, and $\omega$
the K\"ahler form of $S$. Consider the function
$r:\; \Tot B \arrow \R$, $v \arrow |v|^2$.
Using \eqref{_Theta_B_Kah_Equation_},
it is easy to check that $r$ is a K\"ahler potential
on $\Tot B$ (see \cite[(15.19)]{_Besse:Einst_Manifo_}). Let 
\[
X := \{ v \in \Tot B \; | \; |v|^2=1\}
\]
be the circle bundle over $S$ corresponding to $B$.
Denote by $\Tot_0(B)$ the space of non-zero vectors in $B$.
Then, $\Tot_0(B) \cong C(X)$. Since $\Tot_0(B)$ 
is K\"ahler, $X$ is 1-Sasakian.

Considering $X$ as a circle bundle over $S$ and using the Serre's
spectral sequence, we find that $h^1(X) =2$.

Fix $q\in \R^{>1}$.
Let $M= \Tot_0(B) /\sim_q$, where $\sim_q$ is the equivalence
relation generated by $v \sim qv$. Clearly,
$M$ is an LCK manifold. Using the product metric
on $M \cong X \times S^1$, we find that $M$
is actually a Vaisman manifold.

By \eqref{_fibra_over_S^1_cohomo_Equation_},
$h^1(M) = h^1(X)+1 = 3.$ This gives a counterexample
to \ref{_h^1=1_Theorem_} (without the Einstein-Weyl assumption).

The manifold $M$ is equipped with a natural holomorphic projection
$\pi \; M \arrow S$. Clearly, $S$ admits a non-trivial
holomorphic 1-form. Lifting this form to $M$,
we obtain a non-trivial holomorphic form on
a Vaisman manifold. Therefore, for \ref{_holo_forms_vanish_Theorem_},  the 
Einstein-Weyl assumption is also essential.

This example appears (in another language) in 
\cite{_Vaisman:Dedicata_}, under the name of induced
Hopf bundle.

\hfill

\hfill

{\bf Acknowledgements:}
I am grateful to R. Bielawski,  D. Kaledin
and L. Ornea for interesting discussions.

\hfill

% end of small

\end{document}